\pdfoutput=1    

\documentclass[reqno, 12pt]{amsart}
\usepackage{cmap}

\usepackage{amsmath}
\usepackage{amsfonts}
\usepackage{amssymb}
\usepackage{graphicx}

\usepackage{parskip}  

\usepackage{comment}  

\usepackage{url}

\usepackage{ifthen}
\newboolean{BKMRK}
\setboolean{BKMRK}{true}    

\ifthenelse {\boolean{BKMRK}}
  { \usepackage{hyperref} }
  {  }

\ifthenelse {\boolean{BKMRK}}
  { \hypersetup{colorlinks=true, linkcolor=red, citecolor=red, pdfstartview=FitH} }
  {  }

\def\picDblWidth{3.20in}    
\def\figureMinipageWidth{.50\textwidth}  

\theoremstyle{plain}

\numberwithin{equation}{section}

\usepackage[top = .75in, bottom = .75in, left = .75in, right = .75in]{geometry}

\input{glyphtounicode}
\begin{document}

\title{Summing the Curious Series of Kempner and Irwin}
\author{Robert Baillie}


\date{\today}



\begin{abstract}
In 1914, Kempner proved that the series
\[
1/1 + 1/2 + ... + 1/8 + 1/10 + 1/11 + ... + 1/18 + 1/20 + 1/21 + ...
\]
where the denominators are the positive integers that do not contain the digit 9, converges to a sum less than 90.
The actual sum is about 22.92068.
In 1916, Irwin proved, among other things, that the sum of $1/n$ where $n$ has at most a finite number of 9's is also a convergent series.
We show how to compute sums of Irwins' series to high precision. 
For example, the sum of the series
\[
1/9 + 1/19 + 1/29 + 1/39 + 1/49 + ...
\]
where the denominators have exactly one 9, is about 23.04428 70807 47848 31968.
Another example: the sum of $1/n$ where $n$ has exactly 100 zeros is about $10\ln(10) + 1.00745 72171 \cdot 10^{-197} \approx 23.02585$;
note that the first, and largest, term in this series is the tiny 1/googol.
Finally, we discuss a class of related series whose summation algorithm has not yet been developed.
\end{abstract}

\maketitle


\section{Introduction}

The harmonic series
\[
\frac{1}{1} + \frac{1}{2} + \frac{1}{3} + \dots + \frac{1}{n} + \cdots
\]
diverges:
That is, the sum can be made as large as we wish just by adding enough terms.
For example, we can make the sum exceed 10 by adding the first 12367 terms.
We can make the sum exceed 100 by adding 1509 26886 22113 78832 36935 63264 53810 14498 59497 terms \cite{Boas, BoasOEIS, BoasOEIS2}.

But in a 1914 paper titled ``A Curious Convergent Series'', A. J. Kempner \cite{Kempner} proved that the series
\[
\frac{1}{1} + \frac{1}{2} + \frac{1}{3} + \frac{1}{4} + \frac{1}{5} + \frac{1}{6} + \frac{1}{7} + \frac{1}{8} \quad + \quad \frac{1}{10} +
\frac{1}{11} + \frac{1}{12} + \frac{1}{13} + \frac{1}{14} + \frac{1}{15} + \frac{1}{16} + \frac{1}{17} + \frac{1}{18} \quad + \quad \frac{1}{20} + \frac{1}{21} + \cdots
\]
where the denominators are the positive integers \emph{that do not contain the digit 9}, converges to a sum less than 90.
At first glance, this is surprising because Kempner's series looks like the harmonic series with every $10^{\text{th}}$ term deleted.
Kempner's result was indeed curious!

But note that when we reach $1/89$, we delete 11 terms in a row.
Further, in the long run, Kempner's series thins out enough to allow it to converge.
For example, the vast majority of integers having 100 digits contain a 9 somewhere within them,
so most 100-digit numbers do \textit{not} occur as denominators in Kempner's ``no 9'' series.

The author first learned of this result from Hardy and Wright, \cite[pp.\ 120--121]{HardyAndWright}.
This topic has made its way into popular books, for example, by Havil \cite[pp.\ 31--34]{Havil} and Wells \cite[p.\ 81]{Wells},
a ``Mathologer'' YouTube video \cite{Mathologer},
and even into an internet cartoon \cite{smbc-cartoon}.


The sum of Kempner's ``no 9'' series is about 22.92067 66192 64150 34816 \cite{Baillie}.
However, the convergence is so slow that the sum of all $\approx 10^{28}$ terms with denominators less than $10^{30}$
is still less than 22, so this series cannot be accurately summed by brute force addition of terms.

The following generalization also holds: if we delete from the harmonic series all terms whose denominators have any set of numbers with one or more digits, the resulting series also converge.
Schmelzer and Baillie \cite{Schmelzer} showed how to compute sums of these slowly-converging series.
For example, they calculate that the sum of $1/n$ where $n$ has \emph{no even digits}, is about 3.17176 54734 15905 .
They also calculate that the sum of $1/n$ where $n$ has no digit string ``314'', is about 2299.82978 27675 18338 .
Their algorithm also shows that the sum of $1/n$ where $n$ has no ``10'', is about 220.88692 51592 49859 .

In 1916, Frank Irwin \cite{Irwin} generalized Kempner's result in a different way:
he showed that the sum of $1/n$ where $n$ has \emph{at most} any fixed number of occurrences of one or more digits, also converges.
It follows that the sum of $1/n$ where $n$ has \emph{exactly} $n_1$ occurrences of digit $d_1$, $n_2$ occurrences of $d_2$, etc., is a convergent series.
Kempner's ``no 9'' series is simply an Irwin series with one condition placed on the digits, namely: $n_1 = 0, d_1 = 9$.

Irwin's result was also curious.
In fact, Irwin used the same title for his article that Kempner had used two years earlier.

At first glance, it might appear that Irwin's ``one 9'' series, which begins
\[
\frac{1}{9} + \frac{1}{19} + \frac{1}{29} + \cdots
\]
is just the harmonic series with 90\% of its terms deleted.
If that were true, then this series would diverge.
However, if you write down a random 100-digit integer, it will likely have \emph{more than} one 9.
(If the 100 digits are chosen at random, you would expect about \emph{ten} of the digits to be 9's.)

There are $9 \cdot 10^{99}$ integers with 100 digits.
But Equation \eqref{E:Occurrences1} below shows that there are only
\[
9^{100 - 1} + 8 \cdot (100 - 1) \cdot 9^{100 - 2} \approx 2.63 \cdot 10^{96}
\]
integers with 100 digits which have exactly one 9.

Therefore, more than 99.9\% of the integers with 100-digit denominators are excluded from the ``one 9'' series.
In fact, with such a restriction on a single digit, in the long run, the terms thin out enough to make the series converge.

These series studied by Irwin also converge very slowly.
So, a natural question is, ``What are their sums?''
It seems that, until now, no one has accurately computed sums of Irwins' series.
This article shows how to do so.

For example, we will calculate that the sum of $1/n$ where $n$ has exactly one 9 is about 23.04428 70807 47848 31968.

So, we have these amusing results:
\[
\frac{1}{1} + \frac{1}{2} + \frac{1}{3} + \dots + \frac{1}{n} + \cdots
\]
diverges, but
\[
\frac{1}{1} + \dots + \frac{1}{8} \quad + \quad \frac{1}{10} + \dots + \frac{1}{18} \quad + \quad \frac{1}{20} + \dots + \frac{1}{28} \quad + \quad \frac{1}{30} + \dots \approx 22.92068
\]
and
\[
\frac{1}{9} + \frac{1}{19} + \frac{1}{29} + \dots \approx 23.04429.
\]

The sum of these last two series (about 45.96496) is \emph{not} the sum of the harmonic series!
Why?
Observe that $1/99$ is missing from both series.
These last two series, combined, have denominators with \emph{at most} one 9,
while the harmonic series has denominators with arbitrarily many 9's.

Notice that the sum of Irwin's ``one 9'' series is \emph{larger} than the sum of Kempner's ``no 9'' series.
This is surprising because the ``no 9'' series begins with larger terms: $1/1 + 1/2 + \dots + 1/8 + 1/10 + \cdots$,
while the ``one 9'' series begins with smaller terms: $1/9 + 1/19 + 1/29 + 1/39 + \cdots$.

Here's a more elaborate example: the sum of $1/n$ where $n$ has \textit{exactly} one 1, two 2's, three 3's, four 4's, and five 5's
(with no conditions placed on the other digits), is about 0.0046539 02254 05638 15565.
The sum of $1/n$ where $n$ has \textit{at most} one 1, two 2's, three 3's, four 4's and five 5's,
is about 27.56008 29488 96367 05754.

The sum of $1/n$ where $n$ has exactly 100 zeros begins with the miniscule term 1/googol.
However, Table \ref{Ta:s100} in Section ~\ref{S:Googol} shows that, after about $1.2 \cdot 10^{1274}$ terms, the sum of this series exceeds the sum of the ``no 9'' series!

The trick to summing the series of both Kempner and Irwin is to recognize that the denominators in the series obey patterns.
Consider a Kempner sum of $1/n$ where the denominators $n$ have no 5's.
Let $S_i$ be the set of all such $n$ that have $i$ digits.
Then we can form $S_{i+1}$ as follows:
\[
S_{i+1} = \bigcup_{x \in S_i} \{ 10x, 10x+1, 10x+2, 10x+3, 10x+4, 10x+6, 10x+7, 10x+8, 10x+9 \} .
\]
Each $x$ in $S_i$ gives rise to nine numbers in $S_{i+1}$.
Furthermore, the reciprocals of the numbers in $S_{i+1}$ are about $1/10$ as large as those of the corresponding numbers in $S_i$.
Therefore, the sum of reciprocals of elements of $S_{i+1}$ is about $9/10$ as large as the sum of reciprocals of elements of $S_i$.
This argument, made a little more rigorous, explains why Kempner's ``no 9'' series converges.
It also provides a rough bound on the sum of the series.

Section \ref{S:Algorithm} shows how to perform a similar trick with Irwins' series.
 
From now on, a statement like \textbf{``$n$ has $m$ zeros'', means ``$n$ has \emph{exactly} $m$ zeros''.}
These $m$ zeros need not be \emph{consecutive}, that is, together (in a block).
We will also assume that there are no restrictions on any unspecified digits.

Most of the results in this article were calculated with a default to display 15 decimal places.
If the result is displayed to 20 or more decimal places, the digits will be displayed in groups of 5 digits.
All numbers in this article are rounded in the last decimal place unless they are followed by trailing dots.
Finally, the base is assumed to be 10 unless stated otherwise.




\section{What's new in this version?} \label{S:WhatsNew}

The previous version of this article was posted to Math arXiv on March 5, 2023.

This (February, 2024) version has a couple of typos fixed in the text.
No equations were changed.

Section \ref{S:Googol} discusses the sum of $1/n$ where $n$ has 100 zeros.
Some additional information was added to this section.

Section \ref{S:Unsolved} discusses series of the form $1/n$ where $n$ has one or more occurrences of a multi-digit integer.
An incorrect comparison between sums in bases 10 and 100 was deleted.
Some additional information was added to this section, including a rough estimate of the sum of $1/n$ where $n$ has exactly one occurrence of the digit string 35.

The \textit{Mathematica} packages, \verb+irwinSums.m+ and \verb+kempnerSums.m+, were not changed.
These two text files are available as ``ancillary files'' at the arXiv link to this paper at \\
\url{https://arxiv.org/abs/0806.4410} .

The direct links to these two files are: \\
\url{https://arxiv.org/src/0806.4410/anc/irwinSums.m} \\
\url{https://arxiv.org/src/0806.4410/anc/kempnerSums.m} .

\section{The algorithm}\label{S:Algorithm}

Irwins' series, like Kempner's series, converge much too slowly to compute their sums by simply adding terms.
It's easy to add terms having denominators up to, say, 7 digits, but we must use some type of extrapolation procedure to get much beyond that point.
It turns out that we can use sums over denominators with $i$ digits to compute the desired sums over denominators with $i + 1$ digits.
Then we repeat the process.

The algorithm here is derived from the simpler algorithm in \cite{Baillie}, which evaluated the sums of $1/n$ where $n$ has \emph{zero} occurrences of any digit.

For a simple example, consider the set of integers that have exactly one 9.
Among 1-digit numbers, there is just one such number, namely, 9 itself.

To generate the 2-digit numbers with exactly one 9, we perform these two steps:

First, for each 1-digit number $x$ that has \emph{no} 9, we create a 2-digit number with \emph{one} 9 by computing $10 x + 9$.
Here, $x$ has 8 values (1 through 8).
This gives us the eight numbers 19, 29, ..., 79, and 89.

Second, let $x$ be the 1-digit number that has \emph{one} 9.
We compute $10 x + d$, where $d$ is any digit \emph{except} 9.
There are nine such 2-digit values of $10 x + d$, namely, 90, 91, ..., 97, and 98.
These two steps generate all $8 +9 = 17$ 2-digit numbers with exactly one 9.

The same procedure generates 3-digit numbers with one 9:
For each 2-digit $x$ with \emph{no} 9, compute $10x + 9$; there are 72 such $x$.
Then, for each of the 17 2-digit $x$ with \emph{one} 9, compute $10x + d$ where $d$ is any digit \emph{except} 9.
There are $17 \cdot 9 = 153$ such values.
This generates the $72 + 153 = 225$ 3-digit numbers with one 9.

We will now describe the algorithm for when there are \emph{two} conditions for $1/n$ to be included in the series:
namely, where $n$ has exactly $k_1$ occurrences of the digit $d_1$ and exactly $k_2$ occurrences of the digit $d_2$.
The idea extends to conditions on more than two digits, and easily generalizes to bases other than 10.

Define $S(i, k_1, k_2)$ to be the set of positive integers with $i$ digits that have $k_1$ occurrences of digit $d_1$ and $k_2$ occurrences of $d_2$.
We can generate the set $S(i+1, k_1, k_2)$ from $S(i, k_1 - 1, k_2)$, $S(i, k_1, k_2 - 1)$, and $S(i, k_1, k_2)$, in three steps:

(a) For each $x$ in $S(i, k_1 - 1, k_2)$: multiply by 10, then add $d_1$.\\
(b) For each $x$ in $S(i, k_1, k_2 - 1)$: multiply by 10, then add $d_2$.\\
(c) For each $x$ in $S(i, k_1, k_2)$: multiply by 10, then add $d = 0, 1, 2, \dots, 9$, except for $d_1$ and $d_2$.

In the above example with integers having one 9, $d_1 = 9$ and $d_2$ was not present.
We described steps (a) and (c).
Step (b) did not apply.

Step (a) starts with an $i$-digit number having $k_1 - 1$ occurrences of $d_1$ and appends $d_1$ as the final digit, forming an $(i+1)$-digit number having $k_1$ occurrences of $d_1$.
Step (b) does the same for $k_2$ and $d_2$.
In step (c), the $i$-digit numbers already have the desired number of $d_1$ and $d_2$, so in this step, we create $(i+1)$-digit numbers by appending all digits \emph{except} $d_1$ and $d_2$. Together, steps (a), (b), and (c) generate $S(i+1, k_1, k_2)$.
If $k_1$ is 0, we omit step (a).
If $k_2$ is 0, we omit step (b).

\textbf{Example:} Starting with the 2-digit numbers that have one 3 and one 5, we will generate the 3-digit numbers that have one 3 and one 5.
Using the above notation, $i$ = 2, $d_1 = 3$, $d_2 = 5$, and $k_1 = k_2 = 1$.
The three steps are:

(a) $S(2, 0, 1)$ are the 2-digit numbers with no 3, but with one 5. \\  
$S(2, 0, 1) = \{15, 25, 45, 50, 51, 52, 54, 56, 57, 58, 59, 65, 75, 85, 95\}$.\\
For each $x$ in this set, compute $10 x + 3$. \\
We get $\{153, 253, 453, 503, 513, 523, 543, 563, 573, 583, 593, 653, 753, 853, 953\}$.

(b) $S(2, 1, 0)$ are the 2-digit numbers with one 3, but with no 5. \\
$S(2, 1, 0) = \{13, 23, 30, 31, 32, 34, 36, 37, 38, 39, 43, 63, 73, 83, 93\}$.\\
For each $x$ in this set, compute $10 x + 5$. \\
We get $\{135, 235, 305, 315, 325, 345, 365, 375, 385, 395, 435, 635, 735, 835, 935\}$.

(c)  $S(2, 1, 1)$ are the 2-digit numbers with one 3 and one 5. \\
$S(2, 1, 1) = \{35, 53\}$.\\
For each $x$ in this set, compute $10 x + d$, where $d = \{0, 1, 2, 4, 6, 7, 8, 9\}$.\\
We get \{350, 351, 352, 354, 356, 357, 358, 359\} and \{530, 531, 532, 534, 536, 537, 538, 539\}.

The union of these computed sets gives the 46 3-digit numbers having exactly one 3 and one 5.


Section \ref{S:ConfirmingTheCalculations} describes the function \verb+iPartialSumCheck[ ]+.
If you set the $6^{\text{th}}$ parameter to 1, the function prints out the sets
$S(1, k_1, k_2)$, $S(2, k_1, k_2)$, and $S(3, k_1, k_2)$.
This lets you verify that the 46 3-digit numbers we just computed are precisely the 3-digit numbers that belong in the series.

Having computed the elements of the set $S(i, k_1, k_2)$, we now define
\[
t(i, j, k_1, k_2) = \sum_{x\in S(i, k_1, k_2)} \frac{1}{x^j}.
\]
For $j = 1, 2, \dots$, we will show how to compute $t(i+1, j, k_1, k_2)$ by using the values of $t(i, j, k_1 - 1, k_2)$, $t(i, j, k_1, k_2 - 1)$, and $t(i, j, k_1, k_2)$.

If $x$ is an $i$-digit number and $d$ is any digit, then the reciprocal of the $(i+1)$-digit number $10x+d$ can be expanded in powers of $1/x$ :
\[
\frac{1}{10x+d} =  \frac{1}{10x (1+\frac{d}{10x})} = \frac{1}{10x} \sum_{n=0}^{\infty}(-1)^n \left(\frac{d}{10x}\right)^n
= \sum_{n=0}^{\infty}(-1)^n \frac{d^n}{(10x)^{n+1}} \, .
\]

A similar expansion holds for higher powers:
\begin{align}
 \frac{1}{(10x+d)^j} & = \frac{1}{(10x)^j (1+\frac{d}{10x})^j}
 = \frac{1}{(10x)^j} \sum_{n=0}^{\infty}(-1)^n \binom{j+n-1}{n} \left(\frac{d}{10x}\right)^n \notag \\
 & = \frac{1}{(10x)^j} \left( 1 - j \cdot \frac{d}{10x} + \frac{j(j+1)}{2!} \cdot \left(\frac{d}{10x}\right)^2 - \frac{j(j+1)(j+2)}{3!} \cdot \left(\frac{d}{10x}\right)^3 +\ldots \right) \notag \\
 & = \sum_{n=0}^{\infty}(-1)^n \frac{(j+n-1)!}{n! (j-1)!} \frac{d^n}{(10x)^{j+n}} \, . \label{E:negBinomial}
\end{align}

This expansion is the negative binomial series \cite{NegativeBinomialSeries}.
This series converges if $|\frac{d}{10x}| < 1$.

(Note: In Equation \eqref{E:negBinomial}, the first term in the series is 1.
Therefore, if $d = n = 0$, the $0^0$ in the numerator of \eqref{E:negBinomial} is taken to be 1.
This prevents ``indeterminate'' warnings when \emph{Mathematica} tries to evaluate $0^0$.)

Now, recalling step (a), we sum these expansions for all $x$ in $S(i, k_1 - 1, k_2)$.  Call this sum $A$.

\begin{align}
A =& \sum_{x\in S(i, k_1 - 1, k_2)} \frac{1}{(10x + d_1)^j}
= \sum_{x\in S(i, k_1 - 1, k_2)} \sum_{n=0}^{\infty}(-1)^n \frac{(j+n-1)!}{n! (j-1)!} \frac{d_1^n}{(10x)^{j+n}} \label{E:equationA} \\
  &
 = \sum_{n=0}^{\infty}(-1)^n \frac{(j+n-1)!}{n! (j-1)!} \frac{d_1^n}{10^{j+n}}
 \sum_{x\in S(i, k_1 - 1, k_2)} \frac{1}{x^{j+n}} \notag .
\end{align}

We can rearrange the order of this series because it converges absolutely.
To see this, note that the expansion of $1/(10x+d)^j$ but without the $(-1)^n$, is just the expansion of $1/(10x-d)^j$, which also converges.

Moreover, for a given $i$ and $j$, the sum
\[
\sum_{x\in S(i, k_1 - 1, k_2)} \frac{1}{x^{j+n}}
\]
approaches 0 rapidly as $n$ approaches infinity.  In this sum, $x$ ranges over a subset of $i$-digit numbers.  Each value of $x$ is at least $10^{i-1}$, so $1/x \le 1/10^{i-1}$.  So,
\[
\frac{1}{x^{j+n}} \le \frac{1}{ 10^{(i-1)(j+n)} } .
\]
There are fewer than $10^i$ such numbers $x$, so the sum is less than
\[
\frac{10^i}{ 10^{(i-1)(j+n)} } = \frac{1}{\left(10^{i-1}\right)^n} \cdot \frac{10^{i+j}}{10^{ij}} .
\]

Likewise, summing over the sets described in steps (b) and (c), we get
\begin{equation}\label{E:equationB}
B = \sum_{x\in S(i, k_1, k_2 - 1)} \frac{1}{(10x + d_2)^j}
 = \sum_{n=0}^{\infty}(-1)^n \frac{(j+n-1)!}{n! (j-1)!} \frac{d_2^n}{10^{j+n}}
 \sum_{x\in S(i, k_1, k_2 - 1)} \frac{1}{x^{j+n}}
\end{equation}

and

\begin{align}
C = & \sum_{d} \sum_{x\in S(i, k_1, k_2)} \frac{1}{(10x + d)^j} \label{E:equationC} \\
  &
 = \sum_{n=0}^{\infty}(-1)^n \frac{(j+n-1)!}{n! (j-1)!} \cdot \frac{0^n + 1^n + \dots + 9^n - d_1^n - d_2^n}{10^{j+n}}
 \sum_{x\in S(i, k_1, k_2)} \frac{1}{x^{j+n}} \notag .
\end{align}

The sum in \eqref{E:equationC} is over digits $d$ other than $d_1$ or $d_2$.
The factor
\begin{equation*}  
\frac{0^n + 1^n + \dots + 9^n - d_1^n - d_2^n}{10^{j+n}}
\end{equation*}
is less than
\[
\frac{0^n + 1^n + \dots + 9^n}{10^{j+n}}
< \frac{9^n + 9^n + \dots + 9^n}{10^{j+n}}
= \frac{10 \cdot 9^n}{10^{j+n}} = \left( \frac{9}{10} \right)^n \cdot \frac{1}{10^{j-1}}
\]
so for a fixed $j$, this, too, rapidly approaches 0 as $n$ approaches infinity.


Together, the sets generated by steps (a) - (c) form the set $S(i+1, k_1, k_2)$.
Also, $A+B+C = t(i+1, j, k_1, k_2)$.
So, we have just computed a needed sum over $(i+1)$-digit numbers by using sums over $i$-digit numbers.

In order to compute $t(i+1, j, k_1, k_2)$, we used the values of $t(i, j, k_1 - 1, k_2)$ and $t(i, j, k_1, k_2 - 1)$.
But $t(i, j, k_1 - 1, k_2)$, in turn, was computed using $t(i-1, j, k_1 - 2, k_2 - 1)$ and $t(i-1, j, k_1 - 1, k_2 - 2)$.
This means that, for each $i$ and $j$, in order to compute $t(i+1, j, n_1, n_2)$, we must compute all
$(n_1+1)(n_2+1)$ values of $t(i, j, k_1, k_2)$ for $0 \le k_1 \le n_1$ and $0 \le k_2 \le n_2$.

If we set $j = 1$ and add the $t(i, 1, n_1, n_2)$ values over all $i$, then we get the sum of $1/n$ where $n$ has \emph{exactly} $n_1$ occurrences of $d_1$ and \emph{exactly} $n_2$ occurrences of $d_2$.
If we add the $t(i, 1, k_1, k_2)$ values over all $i$, and over all $k_1$ and $k_2$ with $0 \le k_1 \le n_1$ and $0 \le k_2 \le n_2$, then we get the sum of $1/n$ where $n$ has \emph{at most} $n_1$ occurrences of $d_1$ and \emph{at most} $n_2$ occurrences of $d_2$.
This sum may also be of interest, so the function \verb+iSum[ ]+ in the \verb+irwinSums.m+ package prints it out; see Section ~\ref{S:MathematicaImplementation}.

If there is one condition on the digits, say, $n_1$ occurrences of digit $d_1$, then the procedure is similar, except that we omit step (b) above, and ignore all terms involving $d_2$ in equation \eqref{E:equationC}.

If there are three conditions on the digits $(n_1, d_1)$, $(n_2, d_2)$, $(n_3, d_3)$, the procedure is similar, except that we have three steps in place of (a) and (b) above, and that for each $i$ and $j$, we must compute $(n_1+1)(n_2+1)(n_3+1)$ values for the $t$ array.
The time and memory requirements increase accordingly.
Still, even on a personal computer, we can specify a condition for each of the ten digits, provided the product
$(n_1+1)(n_2+1) \dots (n_{10}+1)$ is not too large.
See Examples 3, 4, 5, and 6 in Section ~\ref{S:Examples}.

To summarize, our calculation goes as follows.
We start by calculating $t(i, j, k_1, k_2)$ for $i = 1$, 2, and 3, by explicitly adding the terms whose denominators are in the sets $S(i, k_1, k_2)$, for all $0 \le k_1 \le n_1$ and $0 \le k_2 \le n_2$.
For $i = 1$ and 2, we need only $j = 1$.
For $i = 3$, we explicitly add terms to compute all sums $t(3, j, k_1, k_2)$ for $j \le J$, where $J$ will depend on the number of decimal places desired in the final answer.
Then we begin to extrapolate using equations \eqref{E:equationA} - \eqref{E:equationC}, along with the $t(3, j, k_1, k_2)$ values.
We compute the needed $t(i, j, k_1, k_2)$ for $i$ = 4 and $j = 1,2, \dots, J $.
We then use those values to compute $t(i, j, k_1, k_2)$ for $i = 5$, etc.
We continue until the $t(i, j, k_1, k_2)$ values become small enough to be neglected.
Of course, the more decimal places we want, the larger the range of $i$ and $j$ values we will need to compute.

\section{Mathematica implementation} \label{S:MathematicaImplementation}

\textbf{The irwinSums.m package.}

The text file \verb+irwinSums.m+, a \textit{Mathematica} package that implements the algorithm in Section \ref{S:Algorithm},
can be found at the arXiv link given in Section \ref{S:WhatsNew}.

If you download the file to, say, \verb+C:/math/irwinSums.m+, then you can read that file into \textit{Mathematica} with
\begin{verbatim}
    << C:/math/irwinSums.m
\end{verbatim}

\verb+iSum+ is the function that produced most of the results in this paper.
The general input format for \verb+iSum+ is
\begin{verbatim}
    iSum[ digits, counts for those digits, number of decimal places, base ]
\end{verbatim}
The last two parameters are optional.
If they are omitted, the number of decimals defaults to 15 (or to your current default), and the base defaults to 10.
If the fourth parameter is a base other than 10, the third parameter can be either the number of decimal places, or $-1$ for the current default.
You can change the default number of decimal places with the function \verb+setDefaultDecimals[ ]+.

For example, \verb+iSum[9, 1]+ computes the sum of $1/n$ where $n$ has exactly one 9.
The requested sum is 23.044287080747848.
This is the value that is returned.
The following information is displayed:
\begin{verbatim}
    sum = 23.044287080747848
     sum with exactly 0 occurrences of 9 = 22.920676619264150
     sum with exactly 1 occurrence  of 9 = 23.044287080747848
     sum with at most 1 occurrence  of 9 = 45.964963700011999
\end{verbatim}
As noted in Section \ref{S:Algorithm}, in order to compute the sum of $1/n$ where $n$ has \emph{one} 9,
we also needed to compute the sum where $n$ has \emph{no} occurrence of 9; this result is also displayed.
The last sum, $45.96496 \ldots$, is the sum of $1/n$ where $n$ has \emph{either} no 9's or one 9.


Another example: \verb+iSum[{0, 9}, {2, 3}]+ computes sum of $1/n$ where $n$ has two 0's and three 9's.
The printed output is
\begin{verbatim}
    sum = 3.603022188537689
     sum with at most 2 0's and at most 3 9's = 55.428091573952601
\end{verbatim}
The value returned is the requested sum: $3.60302 \ldots$ .
Because conditions are placed on more than one digit, \verb+iSum+ does not print the sums of all $(2+1)(3+1)=12$ combinations where $n$ has \\
\hspace*{24pt} exactly 0 0's and exactly 0 9's \\
\hspace*{24pt} exactly 0 0's and exactly 1 9 \\
\hspace*{24pt} exactly 0 0's and exactly 2 9's \\
\hspace*{24pt} exactly 0 0's and exactly 3 9's \\
\hspace*{24pt} exactly 1 0 and exactly 0 9's \\
\hspace*{24pt} ... \\
\hspace*{24pt} exactly 2 0's and exactly 3 9's. \\
Only the \emph{total} of these intermediate results, $55.42809 \ldots$, is printed.

There is also a version of \verb+iSum+ which formats the output, breaking digits after the decimal point into groups of 5, for easy reading.
This function is \verb+iSumFormatted+, with the same parameters as \verb+iSum+.
For example, \verb+iSumFormatted[9, 1]+ prints several lines of text, the first of which is, ``sum = 23.04428 70807 47848''.

The easiest way to copy the result is to copy the \emph{printed text}, not the \emph{returned result}.

\verb+iSumFormatted[ ]+ returns a result of \emph{Mathematica} type \verb+NumberForm+.
To copy the \verb+NumberForm+ result:
Highlight the result (23.04428 70807 47848), then use the ``Edit / Copy As / Plain Text'' menu option to copy the result as text.

\textbf{Partial Sums.}
Here are two additional useful functions in the \verb+irwinSums+ package.
The first two and last two parameters are the same as for \verb+iSum[ ]+.
The third parameter in these functions is different.

This function computes the \emph{partial sum} of $1/n$, through denominators $n < 10^p$ :
\begin{verbatim}
    iPartialSum[ digits, counts, p, decimal places, base ]
\end{verbatim}

We saw above that the sum of $1/n$ where $n$ has one 9, is $23.044287080747848$, given by \verb+iSum[9, 1]+.
Let us ask: what is the \emph{partial} sum of terms in this series over denominators less than, say, $10^{18}$
(that is, including all denominators having at most 18 digits)?
This gives us the answer:
\begin{verbatim}
    iPartialSum[9, 1, 18]
\end{verbatim}
This returns $12.870158651317169$, the partial sum of the ``one 9'' series for $n < 10^{18}$.

Suppose you have chosen a threshold, \verb+t+.
In the sum of $1/n$, this function
\begin{verbatim}
    iPartialSumThreshold[ digits, counts, t, decimal places, base ]
\end{verbatim}
tells how many digits are rquired in the denominators to make the partial sum exceed the threshold.


About how far must we go in the series to make the partial sum of the ``one 9'' series reach 20?
We could use trial and error with \verb+iPartialSum[9, 1, k]+ until the sum reaches 20,
But \verb+iPartialSumThreshold[ ]+ eliminates the guesswork:
\begin{verbatim}
    iPartialSumThreshold[9, 1, 20]
\end{verbatim}
This returns a list of four numbers:
\[
\{33, 19.795734087917093, 34, 20.051361976855033\} \, .
\]
These numbers tell us that the partial sum with denominators up to $10^{33}$ is about $19.79573$,
while the partial sum with denominators up to $10^{34}$ is about $20.05136$.
Therefore, we need at lease \emph{some} denominators with 34 digits in order to make the partial sum reach 20.

\bigskip

\textbf{The kempnerSums.m package.}

The \verb+kempnerSums.m+ package implements the algorithm described in \cite{Schmelzer}.
This package computes the sum of $1/n$ where $n$ has no occurrence of some digit or string of digits.
This file can be found at the arXiv link given in Section \ref{S:WhatsNew}.

If you download this text file to, say, \verb+C:/math/kempnerSums.m+, then you can read that file into \textit{Mathematica} with
\begin{verbatim}
    << C:/math/kempnerSums.m
\end{verbatim}

The main function in this package is \verb+kSum[ ]+.

Certain sums can be computed with either the \verb+irwinSums.m+ or the \verb+kempnerSums.m+ package.
For example, after reading in both packages, the sum of $1/n$ where $n$ has no 9 can be computed with either \verb+iSum[9, 0]+ or \verb+kSum[9]+.
Both calculations give 22.920676619264150.

However, each package can compute some things that the other cannot.
For example, \verb+kSum[ ]+ can also compute the sum of $1/n$ where $n$ has no occurrence of a \textit{string} of digits.
For example, the sum of $1/n$ where $n$ has no ``314159'' in base 10 can be computed with \verb+kSum[314159]+.
The result is $2302582.333863782607892$.

On the other hand, to compute the sum of $1/n$ where $n$ has \emph{one} 9,
you must run \verb+iSum[9, 1]+ in the \verb+irwinSums.m+ package.

\section{Confirming the calculations} \label{S:ConfirmingTheCalculations}


The algorithm in Section ~\ref{S:Algorithm} works as follows:
Given conditions on the digits, for example, ``one 3 and one 5'', we directly compute (using ``brute force'') the sum of these $1/n$ where $n$ has 1, 2, and 3 digits;
then, we use equations \eqref{E:equationA} - \eqref{E:equationC} to estimate the sums over $n$ having more than 3 digits.
The result is given by \verb+iSum[{3, 5}, {1, 1}]+, which returns 5.754336150131606.


To begin to verify this result, we can use \verb+iPartialSum[{3, 5}, {1, 1}, 6]+ to display the partial sum for $n < 10^6$.
This function applies Equations \eqref{E:equationA} - \eqref{E:equationC}
to obtain the sum for $10^3 \leq n < 10^6$.
This result is 0.858899223157501.

We can then use \verb+iPartialSumCheck[{3, 5}, {1, 1}, 6]+ to compute directly (i.e., with ``brute force''),
the sum over those $n < 10^6$ that have one 3 and one 5.
This result is 0.858899223157501.
This function works by examining the digits of each $n < 10^6$, and adding $1/n$ for only those $n$ that have exactly one 3 and one 5.
So, this function is independent of the algorithm described in Section \ref{S:Algorithm}.
However, \verb+iPartialSumCheck[ ]+ is also much slower.

\verb+iPartialSumCheck[{3, 5}, {1, 1}, 6]+ prints results up to $n < 10^k$ for \emph{each} $k$, $1 \leq k \leq 6$.
These values agree with the six corresponding values of \verb+iPartialSum[{3, 5}, {1, 1}, k]+.

Recall the example in Section \ref{S:Algorithm}, where we used the two 2-digit numbers that have one 3 and one 5
to generate the 46 3-digit numbers that have one 3 and one 5.
If you want to verify that we obtained \emph{exactly} the \emph{correct} values of $n < 10^3$, you can
print out the list of all such $n < 10^3$ by setting the $6^{\text{th}}$ parameter to 1:
\begin{verbatim}
  iPartialSumCheck[{3, 5}, {1, 1}, 3, -1, 10, 1]
\end{verbatim}
(The list is printed only for $n < 10^3$, so there's no need for the $3^{\text{rd}}$ parameter to be more than 3.)

In Section \ref{S:Algorithm}, we also computed the 2- and 3-digit numbers having one 9.
We can verify these 242  2- and 3-digit values of $n$ with \verb+iPartialSumCheck[9, 1, 3, -1, 10, 1]+.

Most of the series in this paper converge very slowly.
However, in a few special cases, they converge rapidly enough that we can compute their sums directly.
These cases can serve as a check on the algorithm.
For example, the sum of $1/n$ where $n$ has \emph{only} 1's in base 2 is
\[
\frac{1}{1} + \frac{1}{3} + \frac{1}{7} + \frac{1}{15} +  \dots = \sum_{m=1}^{\infty} \frac{1}{2^m-1}
\]
which converges rapidly to about 1.60669 51524 15291 76378.
(This constant is known as the Erd{\H o}s-Borwein constant \cite{EB-Constant}.
By a theorem of Borwein, \cite[Theorem 1]{Borwein}, this number is irrational;
see Section \ref{S:Irrat}.)

\emph{Mathematica} gives an \emph{exact} value for this sum, along with a numerical approximation.
These are
\[
\frac{\log (2)-\psi _{\frac{1}{2}}(1)}{\log (2)} \approx 1.60669 \ 51524 \ 15291 \ 76378 \, ,
\]
where $\psi_q(z)$ is the $q$-digamma function (\verb+QPolyGamma[z, q]+ in \emph{Mathematica}).
This exact expression comes from \cite{Borwein-Borwein} and \cite{Weisstein-Lambert}.

This series has no 0's in base 2, so the sum can be calculated to 20 decimals with
\begin{verbatim}
    iSum[0, 0, 20, 2]
\end{verbatim}
which returns 1.60669 51524 15291 76378.
Because we are specifying a value of 2 in the fourth parameter, we cannot omit the third parameter, the number of decimal places (20).
(To use the current \emph{default} value for the number of decimal places, you can set the third parameter to $-1$.)

The sum of $1/n$ where $n$ has \emph{one} 0 in base 2 is (to 25 decimals) \verb+iSum[0, 1, 25, 2]+ = 1.46259 07350 44364 69954 61454.
Writing out the denominators in binary, one can see that the sum is the rapidly-converging series
\[
\sum_{m=2}^{\infty} \sum_{k=0}^{m-2} \frac{1}{2^m-1-2^k} \approx 1.46259 \ 07350 \ 44364 \ 69954 \ 61454 \, .
\]
For example, if $m = 5$, $2^5  - 1$ has five 1's in base 2;
the four numbers having one 0 are: $30 = 2^5 - 1 - 2^0 = 11110_2$, $29 = 2^5 - 1 - 2^1 = 11101_2$, $27 = 2^5 - 1 - 2^2 = 11011_2$, and $23 = 2^5 - 1 - 2^3 = 10111_2$.

\verb+iPartialSumCheck[..., ..., 6]+ is slow because it examines all numbers up through 6 digits, that is, up to $10^6$.
But if you use a smaller base, you can quickly check results that have more digits.
$5^8$ is only 390625, so if you use base 5, you can quickly check sums up through 8 digits.

For example, the sum of $1/n$ where $n$ has two 1's in base 5, is computed with
\verb+iSum[1, 2, -1, 5]+.
The sum is 8.058128592786750.
\verb+iPartialSum[1, 2, 8, -1, 5]+ computes the sum through 8 base-5 digits.
This sum is  2.719119119003371.
Both of these calculations used the algorithm in Section \ref{S:Algorithm}.
\verb+iPartialSumCheck[1, 2, 8, -1, 5]+ computes sums through 8 base-5 digits \emph{without} using the algorithm.
This result is 2.71911911900337.

Finally, when we calculate the sum of $1/n$ where $n$ has one occurrence of a digit $d$, as an intermediate step, the algorithm also computes the sum of $1/n$ where $n$ has zero occurrences of $d$.
When we do these calculations to 5000 decimals for each of the ten digits, the sums for zero occurrences match the corresponding 5000 decimal place sums obtained by the algorithm in \cite{Schmelzer}.

\section{Counting terms in the series} \label{S:Counting}

We will derive formulas that count the number of terms in series like the sum of $1/n$ where $n$ has
$m$ occurrences of any digit.
Our formulas will vary depending on whether the digit in question is 0.
This is because any digit \emph{except} 0 can be the leading digit of a positive integer.
However, the counts will be the same for all non-zero digits.
To make the argument more concrete, we will work through examples using the digit 9,
but those formulas will also apply to all digits $1 \leq d \leq 9$.

Let $k$ be the number of digits in a positive integer $n$.

Also, let $K_d(k, m)$ be the count of integers with $k$ digits, having $m$ occurrences of the digit $d$.

\textbf{Zero occurrences of a digit.}

We can describe this situation as, `$n$ is missing some digit $d$'.

Among $k$-digit numbers having no 0's, there are 9 choices (1 through 9) for the leading digit.
For each of these choices, there are also 9 choices (1 through 9) for each of the other $k - 1$ digits.
So, the number of $n$ with $k$ digits, but with no 0's, is $9^k$.

If a $k$-digit number has no 9's, then there are 8 choices (1 through 8) for the leading digit.
For each of these choices, each of the other $k - 1$ digits can take 9 values (0 through 8).
Therefore, there are $8 \cdot 9^{k - 1}$ numbers with $k$ digits having no 9's.

So, the number of $n$ which have $k$ digits, but which have zero occurrences of some digit, is
\begin{equation} \label{E:Occurrences0}
K_d(k, 0) = 8 \cdot 9^{k - 1} \ \text{ if } d \neq 0, \ \text{ and } \ K_0(k, 0) = 9^k \, .
\end{equation}

\textbf{Zero occurrences of $j \geq 2$ different digits.}

We will work through examples with $j = 2$.
There are two cases: (1) one of the missing digits is 0, and (2) none of the missing digits is zero.
So, we'll count (1) the $n$ with $k$ digits which have no 0 and no 9, and
(2) the $n$ with $k$ digits which have no 8 and no 9.

Case 1: $n$ has no 0 and no 9.
There are 8 choices (1 through 8) for each of the $k$ digits.
Therefore, the number of such $n$ is $8^k$.

Case 2: $n$ has no 8 and no 9.
There are 7 choices (1 through 7) for the leading digit.
For each of these choices, there are 8 choices (0 through 7) for each of the remaining $k - 1$ digits.
Therefore, the number of such $n$ is $7 \cdot 8^{k - 1}$.


It is easy to generalize this argument to the case where $n$ is missing $j$ different digits.

The number of $n$ with $k$ digits, but which have zero occurrences of $j$ distinct digits, is
\begin{align} \label{E:mMissingDigits}
(10 - j)^k & \medspace \text{ if one of the missing digits is 0, or } \notag \\
(10 - j - 1) \cdot (10 - j)^{k - 1} & \medspace \text{ if none of the missing digits is zero} \, .
\end{align}
Note that in the special case $j = 1$, this is equivalent to Equation \eqref{E:Occurrences0}.

\textbf{One occurrence of a digit.}

For numbers that have one occurrence of a non-zero digit, say, 9, there are two cases.

Case 1: The leading digit is 9.
\emph{Each} of the other $k - 1$ digits can be 0 through 8 (9 values), so the number of these $n$ is $9^{k - 1}$.

Case 2: The leading digit is \emph{not} 9.
There are 8 choices (1 through 8) for the leading digit.
Among the remaining $k - 1$ digits, there is one 9.
For each of the 8 choices of the leading digit:
\begin{itemize}
 \setlength\itemsep{-0.5em}
 \item if the 2nd digit is 9, there are $9^{k - 2}$ ways to fill the $k - 2$ digits 3 through $k$
 \item if the 3rd digit is 9, there are $9^{k - 2}$ ways to fill the $k - 2$ digits 2, and 4 through $k$
 \item . . .
 \item if the $k$th digit is 9, there are $9^{k - 2}$ ways to fill the $k - 2$ digits 2 through $k - 1$.
\end{itemize}
In Case 2, there are $\binom{k-1}{1} = k - 1$ possible locations for the 9, so there are $8 \cdot (k - 1) \cdot 9^{k - 2}$ possible values for $n$.

Therefore, the total number of $n$ with $k$ digits, one of which is 9, is:
\[
9^{k - 1} + 8 \cdot (k - 1) \cdot 9^{k - 2} \, .
\]

For numbers that have one occurrence of 0, Case 1 does not apply.
In Case 2, there are 9 choices (1 through 9) for the leading digit.
The rest of the argument is the same, so the number of such $n$ is
\[
9 \cdot (k - 1) \cdot 9^{k - 2} = (k - 1) \cdot 9^{k - 1} \, .
\]

Combining these results, we have
\begin{equation} \label{E:Occurrences1}
K_d(k, 1) = 9^{k - 1} + 8 \cdot (k - 1) \cdot 9^{k - 2} \ \text{ if } d \neq 0, \ \text{ and } \
K_0(k, 1) = (k - 1) \cdot 9^{k - 1} \, .
\end{equation}

\textbf{Two occurrences of a digit.}

If there are two occurrences of a nonzero digit, say, two 9's, we have two cases.

Case 1: The leading digit is 9.
Among the remaining $k - 1$ digits, there is one 9.
There are $\binom{k-1}{1} = k  - 1$ ways to choose which digit position is 9.
Each of the remaining $k - 2$ digits can have 9 values (0 through 8).
The number of these $n$ is therefore $(k - 1) \cdot 9^{k - 2}$.

Case 2: The leading digit is \emph{not} 9.
The leading digit can have 8 values (1 through 8).
For each of these 8 choices of the leading digit, we can select the digit positions of the two 9's in $\binom{k-1}{2}$ ways.
Each of the remaining $k - 3$ digits can be have 9 values (0 through 8).
The number of these $n$ is $8 \cdot \binom{k-1}{2} \cdot 9^{k - 3}$.

If, instead, there are two 0's, only Case 2 applies.
For each of the 9 possible values of the leading digit, we can choose the two digit locations of the 0's in $\binom{k-1}{2}$ ways.
Each of the remaining $k - 3$ digits can have 9 values.

Combining these results, we have
\begin{equation*} 
K_d(k, 2) = (k - 1) \cdot 9^{k - 2} + 8 \cdot \binom{k - 1}{2} \cdot 9^{k - 3} \ \text{ if } d \neq 0, \ \text{ and } \
K_0(k, 2) = 9 \cdot \binom{k - 1}{2} \cdot 9^{k - 3} \, .
\end{equation*}

\textbf{$m$ occurrences of a digit.}



These arguments generalize.
With $m$ 9's, if the leading digit \emph{is} 9, then there are $\binom{k - 1}{m - 1}$ positions for the remaining $(m - 1)$ 9's.
If the leading digit is \emph{not} 9, then there are $\binom{k - 1}{m}$ positions for the $m$ 9's.
The order of these positions does not matter, because we are choosing locations for all copies of the \emph{same} digit.
There are 9 choices (0 through 8) for every digit that is \emph{not} a nine.

So, for any non-zero digit $d$:
\begin{equation} \label{E:OccurrencesNotZeroM}
K_d(k, m) = \binom{k - 1}{m - 1} \cdot 9^{k - m} + 8 \cdot \binom{k - 1}{m} \cdot 9^{k - m - 1} \, ,
\end{equation}
where $\binom{k - 1}{m - 1} = 0$ if $m = 0$.
For $m$ occurrences of the digit 0,
\begin{equation} \label{E:OccurrencesZeroM}
K_0(k, m) = 9 \cdot \binom{k - 1}{m} \cdot 9^{k - m - 1} \, .
\end{equation}
These equations hold for all $m \geq 0$, so they generalize the previous expressions in this Section.

\textbf{One occurrence each of two non-zero digits.}

Here, we will count the number of $n$ which have $k$ digits, and which have one 8 and one 9.
As before, this reasoning holds for any two distinct non-zero digits.

There are \emph{three} cases.

Case 1: The leading digit is 8.
Among the remaining $k - 1$ digits, one is 9, while each of the other $k - 2$ digits has 8 possible values (0 through 7).
There are $\binom{k-1}{1} = k - 1$ ways to choose the 9's digit position.
The remaining $k - 2$ digits have $8^{k - 2}$ possible values.
Therefore, there are $\binom{k-1}{1} \cdot 8^{k - 2}$ such $n$.

Case 2: The leading digit is 9.
Interchanging the roles of the 8 and 9, there are \emph{another} $\binom{k-1}{1} \cdot 8^{k - 2}$ such $n$ with leading digit 9.

Case 3: The leading digit is neither 8 nor 9 (and, as always, is not 0).
The leading digit has 7 possible values (1 through 7).
Among the remaining $k - 1$ digits, there are $\binom{k - 1}{2}$ ways to choose the locations for the 8 and the 9.
But here, \emph{the order of the 8 and 9 does matter}.
There are $2!$ possible orders for the two digits 8 and 9.
Finally, the remaining $k - 3$ digits each have 8 possible values (from 0 through 7).


Therefore, the total number of $n$ having $k$ digits, and which have exactly one each of two non-zero digits, is
\begin{equation} \label{E:OccurrencesOneAOneB}
2 \binom{k-1}{1} \cdot 8^{k - 2} + 7 \cdot 2! \cdot \binom{k - 1}{2} \cdot 8^{k - 3} \, .
\end{equation}

\textbf{One occurrence each of three non-zero digits.}

Finally, we will state the following without proof.
The derivation is similar to that of Equation \eqref{E:OccurrencesOneAOneB}.
The number of $n$ having $k$ digits, and which have exactly one each of three non-zero digits, is
\begin{equation} \label{E:OccurrencesOneAOneBOneC}
3 \cdot 2! \cdot \binom{k - 1}{2} \cdot 7^{k - 3} + 6 \cdot 3! \cdot \binom{k - 1}{3} \cdot 7^{k - 4} \, .
\end{equation}

\textbf{Checking Our Results.}
We can confirm Equations \eqref{E:Occurrences0} - \eqref{E:OccurrencesOneAOneBOneC} as follows.
For example, consider the sum of $1/n$ where $n$ has one 9.
\verb+iPartialSumCheck[9, 1, 7]+ computes partial sums to $10^k$, for $1 \leq k \leq 7$.
This function also prints out the \emph{counts} of these $n$ that have 1, 2, 3, 4, 5, 6, and 7 digits;
these counts are: 1, 17, 225, 2673, 29889, 321489, and 3365793.
These counts match the values given by Equations \eqref{E:Occurrences1} and \eqref{E:OccurrencesNotZeroM}, for $1 \leq k \leq 7$.

Likewise, Equation \eqref{E:OccurrencesOneAOneBOneC} can be checked with \verb+iPartialSumCheck[{7, 8, 9}, {1, 1, 1}, 7]+.

\textbf{Generalizations.}
Equations \eqref{E:Occurrences0} - \eqref{E:OccurrencesOneAOneBOneC} generalize to other bases.
For example, for base $b$, you would replace the 9's with $b - 1$, the 8's with $b - 2$, the 7's with $b - 3$, and so on.

\textbf{Counts for Various Series.}

\begin{table}[ht]
 \begin{center}
  \begin{tabular}{ l r r r r r r }
  $k$   &  no 0 & no 9 &   one 0  & one 9 & two 0 & two 9 \\ \hline
 1    &           9 &          8 &          0 &          1 &        0 &            0 \\
 2    &          81 &         72 &          9 &         17 &        0 &            1 \\
 3    &         729 &        648 &        162 &        225 &        9 &           26 \\
 4    &        6561 &       5832 &       2187 &       2673 &       243 &         459 \\
 5    &       59049 &      52488 &      26244 &      29889 &      4374 &        6804 \\
 6    &      531441 &     472392 &     295245 &     321489 &     65610 &       91125 \\
 7    &     4782969 &    4251528 &    3188646 &    3365793 &    885735 &     1141614 \\
 8    &    43046721 &   38263752 &   33480783 &   34543665 &   11160261 &   13640319 \\
 9    &   387420489 &  344373768 &  344373768 &  349156737 &  133923132 &  157306536 \\
10    &  3486784401 & 3099363912 & 3486784401 & 3486784401 & 1549681956 & 1764915561 \\
15    & $2.059 \cdot 10^{14}$ & $1.830 \cdot 10^{14}$ & $3.203 \cdot 10^{14}$ & $3.076 \cdot 10^{14}$ & $2.313 \cdot 10^{14}$ & $2.412 \cdot 10^{14}$ \\
20    & $1.216 \cdot 10^{19}$ & $1.081 \cdot 10^{19}$ & $2.567 \cdot 10^{19}$ & $2.417 \cdot 10^{19}$ & $2.567 \cdot 10^{19}$ & $2.567 \cdot 10^{19}$ \\
  \end{tabular}
   \caption{Counts of $n$ having exactly $k$ digits, for various series}
   \label{Ta:Counts}
  \end{center}
\end{table}

Table \ref{Ta:Counts} shows, for various series (no 0, no 9, one 0 one 9, etc.), how many $n$ there are with $k$ digits.

Notice that, in the row for $k = 10$, three numbers are the same.
If we compute the \emph{exact} values for $k = 20$ and $k = 30$, we find, using the notation of Equations \eqref{E:OccurrencesNotZeroM} and \eqref{E:OccurrencesZeroM}, that:
\begin{itemize}
\setlength\itemsep{-0.5em}
 \item for $k = 10$:  $K_0(10, 0) = K_0(10, 1) = K_9(10, 1)$
 \item for $k = 20$:  $K_0(20, 1) = K_0(20, 2) = K_9(20, 2)$
 \item for $k = 30$:  $K_0(30, 2) = K_0(30, 3) = K_9(30, 3)$.
\end{itemize}

The general claim would be that, for any integer $i \geq 1$,
\[
K_0(10 i, i - 1) = K_0(10 i, i) = K_9(10 i, i) \, .
\]

For $k = 10$, if we write out the expressions for $K_0(10, 0)$, $K_0(10, 1)$, and $K_9(10, 1)$, we claim that
\[
9 \cdot \binom{10 - 1}{0} \cdot 9^{10 - 0 - 1}
= 9 \cdot \binom{10 - 1}{1} \cdot 9^{10 - 1 - 1}
= \binom{10 - 1}{1 - 1} \cdot 9^{10 - 1} + 8 \cdot \binom{10 - 1}{1} \cdot 9^{10 - 1 - 1} \, .
\]

The left-hand side is $9 \cdot (1) \cdot 9^{10 - 1} = 9^{10}$.
The middle part is $9 \cdot (9) \cdot 9^{10 - 2} = 9^{10}$.
The right-hand side is $(1) \cdot 9^{10 - 1} + 8 \cdot (9) \cdot 9^{10 - 2} 
= 9^9 + 8 \cdot 9^9= 9^{10}$.

The interested reader is welcomed to verify the other claims.

\begin{table}[ht]
 \begin{center}
  \begin{tabular}{ l r r r }
  $k$   &  one 9 & one 8, one 9 &  one 7, one 8, one 9 \\ \hline
 1    &            1  &           0  &          0  \\
 2    &           17  &           2  &          0  \\
 3    &          225  &          46  &          6  \\
 4    &         2673  &         720  &        162  \\
 5    &        29889  &        9472  &       2772  \\
 6    &       321489  &      112640  &      38220  \\
 7    &      3365793  &     1253376  &     463050  \\
 8    &     34543665  &    13303808  &    5142942  \\
 9    &    349156737  &   136314880  &   53647944  \\
10    &   3486784401  &  1358954496  &  533655864  \\
15    &  $3.076 \cdot 10^{14}$ & $1.029 \cdot 10^{14}$ & $3.347 \cdot 10^{13}$ \\
20    &  $2.417 \cdot 10^{19}$ & $6.075 \cdot 10^{18}$ & $1.398 \cdot 10^{18}$ \\
50    &  $2.551 \cdot 10^{48}$ & $4.808 \cdot 10^{46}$ & $5.338 \cdot 10^{44}$ \\
100    &  $2.627 \cdot 10^{96}$ & $2.765 \cdot 10^{92}$ & $7.881 \cdot 10^{87}$ \\
  \end{tabular}
   \caption{Counts of $n$ having exactly $k$ digits, for three series}
   \label{Ta:Counts2}
  \end{center}
\end{table}

Table \ref{Ta:Counts2} shows counts for three series, which have a single occurrence of one, two distinct, and three distinct non-zero digits.
Again, these \emph{counts} will be the same regardless of \emph{which} distinct non-zero digits these are.
However, the \emph{sums} do depend on which digits are involved.
The sums of $1/n$ where $n$ has (a) one 9, (b),  one 8 and one 9, and (c) one 7, one 8, and one 9, are about
23.044287080747848, 5.756763686415333, and 1.705545235515835, respectively.

\section{Examples} \label{S:Examples}

\textbf{Example 1 (a).}
Consider the sum of $1/n$ where $n$ has no 9's.
We can compute this sum with \verb+iSum[9, 0]+.
The result is 22.920676619264150.

How far must we go in the series until the partial sum reaches 22?
The answer is given by
\begin{verbatim}
  iPartialSumThreshold[9, 0, 22]
\end{verbatim}
This returns approximately \{30, 21.971, 31, 22.066\}.
These numbers mean that the sum over $n < 10^{30}$ (where $n$ has no 9's) is about 21.971,
but the sum over $n < 10^{31}$ is about 22.066.

About how many \emph{terms} are needed to make the sum exceed 22?

Equation \eqref{E:OccurrencesNotZeroM} with $m = 0$ shows that there are $8 \cdot 9^{k-1}$ numbers with $k$ digits which have no 9's.
So, the total numbers of such $n$ having 1 through 30 or 31 digits, are
\[
\sum_{k = 1}^{30} 8 \cdot 9^{k-1} \approx 4.239 \cdot 10^{28}
 \, , \quad \text{and} \quad
\sum_{k = 1}^{31} 8 \cdot 9^{k-1} \approx 3.815 \cdot 10^{29} \, .
\]
Therefore, the sum of $4 \cdot 10^{28}$ terms is a little less than 22,
but the sum of $4 \cdot 10^{29}$ terms is a little more than 22.

\textbf{Example 1 (b).}
Consider the sum of $1/n$ where $n$ has one 9.
We can compute this sum with \verb+iSum[9, 1]+.
The result is $23.044287080747848$.

How far must we go before the partial sum reaches 22?
\verb+iPartialSumThreshold[9, 1, 22]+ returns approximately $\{46, 21.964, 47, 22.055\}$.
So, the sum over $n < 10^{46}$ is about $21.964$, while the sum over $n < 10^{47}$ is about $22.055$.

About how many terms are there in this series for $n < 10^{46}$ or $n < 10^{47}$?

Equation \eqref{E:OccurrencesNotZeroM} tells how many $k$-digit numbers have one 9.
With $m = 1$, this Equation reduces to
\[
9^{k - 1} + 8 \cdot (k - 1) \cdot 9^{k - 2} \, .
\]
We sum this for $k = 1$ through $k = 46$ and $k = 47$.
We conclude that there are about $4.015 \cdot 10^{44}$ terms in this series for $n < 10^{46}$,
and about $3.692 \cdot 10^{45}$ terms for $n < 10^{47}$.

Therefore, in the sum of $1/n$ where $n$ has exactly one 9,
the sum of the first $4.015 \cdot 10^{44}$ terms with $n < 10^{46}$ is a little less than 22,
while the sum of the first $3.692 \cdot 10^{45}$ terms with $n < 10^{47}$ is a little more than 22.

It is interesting to compare this to Example 1 (a).
There, the sum of the ``no 9'' series reached 22 after ``only'' about $4 \cdot 10^{29}$ terms.

\textbf{Example 2.}
The Cantor ``ternary'' (or ``middle third'') set is obtained recursively as follows.
Start with the interval $[0, 1]$.
Remove the middle third, that is, the open interval $(1/3, 2/3)$.
Then, from the remaining intervals, $[0, 1/3]$ and $[2/3, 1]$ remove the middle thirds, that is, $(1/9, 2/9)$ and $(7/9, 8/9)$.
Continue this process recursively, removing the open middle third from each of the remaining intervals.
What remains is the Cantor set.
The numbers in $[0, 1]$ that remain are precisely those with no 1 in their base 3 representation.

The set of positive \emph{integers} that have no 1 in their base 3 representation begins
$2 = 2_3, 6 = 20_3, 8 = 22_3, 18 = 200_3, 20 = 202_3, 24 = 220_3, 26 = 222_3$.
The sum of $1/n$ where $n$ has no 1's in base 3 is given by \verb+iSum[1, 0, -1, 3]+ $= 1.341426555483088$.

\textbf{Example 3 (a).}
Consider the sum of $1/n$ where $n$ has exactly one 8 and one 9.
For brevity, we'll use 10 decimal places in this example.
\begin{verbatim}
  iSum[{8, 9}, {1, 1}, 10]
\end{verbatim}
The printed output is
\begin{verbatim}
sum = 5.7567636864
 sum with at most 1 8 and at most 1 9 = 28.5987554834
\end{verbatim}
Where does that 28.5987554834 come from?

28.5987554834 is the sum four series, where $n$ has (a) no 8 and no 9, (b) exactly one 8 and no 9,
(c) no 8 and exactly one 9, and (d) exactly one 8 and exactly one 9.

In other words, it is the sum of the following four series:
\begin{verbatim}
  iSum[{8, 9}, {0, 0}, 10]
  iSum[{8, 9}, {1, 0}, 10]
  iSum[{8, 9}, {0, 1}, 10]
  iSum[{8, 9}, {1, 1}, 10]
\end{verbatim}
The sums of these four series are, respectively, 11.2915816168, 5.8240865454, 5.7263236348, and 5.7567636864 .
The sum of these numbers is 28.5987554834.
(The default `print level' is 1.
But if you use \verb+setPrintLevel[2]+ to increase it to 2, all of these these sums will be printed out.)

In general, \verb+iSum[{8, 9}, {a, b}]+ would print out a line similar to
\begin{verbatim}
 sum with at most a 8's and at most b 9's = x
\end{verbatim}
In this case, $x$ would be the sum of $(a + 1)(b + 1)$ different series, which would consist of every combination of between 0 and $a$ occurrences of 8
and between 0 and $b$ occurrences of 9.

\textbf{Example 3 (b).}
In the introduction, we claimed that sum of $1/n$ where $n$ has exactly one 1, two 2's, three 3's, four 4's, and five 5's, is about 0.0046539 02254 05638 15565,
and that the sum of $1/n$ where $n$ has \textit{at most} one 1, two 2's, three 3's, four 4's and five 5's,
is about 27.56008 29488 96367 05754.

If you run
\begin{verbatim}
  iSum[{1, 2, 3, 4, 5}, {1, 2, 3, 4, 5}, 20]
\end{verbatim}
The \emph{printed} output is
\begin{verbatim}
  sum = 0.0046539022540563815565
   sum with at most 1 1, ..., and at most 5 5's = 27.56008294889636705754
\end{verbatim}

The \emph{number of series} whose sums add up to 27.56008 ... is $(1+1)(2+1)(3+1)(4+1)(5+1) = 720$.

\textbf{Example 4 (a).}
Consider the sum of $1/n$ where $n$ has exactly \textit{one} of each digit.
Because we are limiting the number of occurrences of every digit, this series will have only a finite number of terms.
In every denominator of this series, each digit occurs exactly once, so every denominator has exactly 10 digits, all distinct.
There are $10! \times 9/10 = 3265920$ such numbers.
The approximate sum of this finite series can be calculated to 20 decimals with
\begin{verbatim}
  iSum[{0, 1, 2, 3, 4, 5, 6, 7, 8, 9}, {1, 1, 1, 1, 1, 1, 1, 1, 1, 1}, 20].
\end{verbatim}
The printed output is
\begin{verbatim}
    sum = .00082589034791925293861
     sum with at most 1 0, ..., and at most 1 9 = 8.92994817475544342417
\end{verbatim}

\textbf{Example 4 (b).}
In order obtain the sum in Example 4 (a), our algorithm first computes the $2^{10}$ sums over those $n$ which have all combinations of either zero or one occurrence of each of the ten digits.
Together, these $n$ comprise the positive integers that have distinct digits.
There are 8877690 of these integers from 1 through 9876543210.
Elementary problem E2533 in the \textit{American Mathematical Monthly} \cite{Pondiczery} asks for the sum of their reciprocals.
Example 4 (a) shows that this sum is about 8.92995, a more precise answer than the one given in \cite{Foregger}.
Interestingly, Eric Weisstein \cite{Weisstein-Digit} computed the \textit{exact} value of this finite sum, a fraction whose numerator and denominator have 14816583 and 14816582 digits, respectively.

Table \ref{Ta:OneOfEachDigit} shows a similar calculation for bases 2 through 10.
For example, in base 3, there are four $n$ that have exactly one 0, one 1, and one 2.
These are $102_3 = 11$, $120_3 = 15$, $201_3 = 19$, and $210_3 = 21$.
The sum of their reciprocals is $1886/7315 \approx 0.25782 63841$.
Our algorithm gives \verb+iSum[{0, 1, 2}, {1, 1, 1}, 20, 3] = 0.25782638414217361586+ .

\begin{table}[ht]
 \begin{center}
  \begin{tabular}{ r l l }
  $b$ &     \textit{exactly} one of each digit  & \textit{at most} one of each digit \\ \hline
   2  &    0.5    &  1.5  \\
   3  &    0.25782 63841 42173 61586  &  2.60068 35269 99316 47300  \\
	 4  &    0.12978 48084 06223 61677  &  3.60808 68029 86413 68569  \\
	 5  &    0.061089 50418 58837 22653  &  4.56754 45346 30053 32968  \\
	 6  &    0.027277 91550 93573 00827  &  5.49201 22502 32697 89881  \\
	 7  &    0.011747 43887 45155 38097  &  6.38730 21972 13363 81858  \\
	 8  &    0.0049323 10351 55588 69283  &  7.25702 34020 72060 68918  \\
	 9  &    0.0020326 81860 44365 59280  &  8.10385 28180 85605 31093  \\
	10  &    0.00082589 03479 19252 93861  &  8.92994 81747 55443 42417  \\
  \end{tabular}
   \caption{The sum of $1/n$ where $n$ has \textit{exactly}, or \textit{at most}, one of each digit, base $b$}
   \label{Ta:OneOfEachDigit}
  \end{center}
\end{table}

\textbf{Example 5 (a).}
Consider the sum of $1/n$ where $n$ has exactly \textit{two} of each digit.
Here, all denominators in this finite series have 20 digits.
This sum can be computed to 20 decimals with
\begin{verbatim}
  iSum[{0, 1, 2, 3, 4, 5, 6, 7, 8, 9}, {2, 2, 2, 2, 2, 2, 2, 2, 2, 2}, 20]
\end{verbatim}
The result is 0.000054406 21942 90990 91465 .

\textbf{Example 5 (b).}
As part of the calculation in 5 (a), the algorithm also computed the $3^{10}$ sums of $1/n$ where $n$ has zero, one, or two occurrences of each of the ten digits.
Together, these sums comprise a finite sum that terminates after 20-digit denominators.
The algorithm prints the total of these sums, that is, the sum of $1/n$ where $n$ has at most two of every digit.
This sum is about 20.58988 67749 18085 64961.

Table \ref{Ta:TwoOfEachDigit} shows the results of similar calculations for bases 2 through 10.
For base 2, the two sums are $1/9 + 1/10 + 1/12 = 53/180$ and $1/1 + 1/2 + 1/3 + 1/4 + 1/5 + 1/6 + 1/9 + 1/10 + 1/12 = 247/90$, respectively.

\begin{table}[ht]
 \begin{center}
  \begin{tabular}{ r l l }
  $b$ &     \textit{exactly} two of each digit  & \textit{at most} two of each digit \\ \hline
   2  &    0.29444 44444 44444 44444    &  2.74444 44444 44444 44444  \\
   3  &    0.13728 42813 53796 22429  &  5.03830 16073 84291 30946  \\
	 4  &    0.053227 67916 73449 82805  &  7.31560 12470 91568 19264  \\
	 5  &    0.018595 61585 71205 43656  &  9.58154 13870 01674 99146  \\
	 6  &    0.0061233 38585 38739 23485  &  11.82888 19161 80409 98185  \\
	 7  &    0.0019419 37342 50659 60903  &  14.05401 73172 65876 98270  \\
	 8  &    0.00060011 95484 45774 43387  &  16.25577 41877 60599 32348  \\
	 9  &    0.00018199 15685 46436 33810  &  18.43419 29090 32510 52124  \\
	10  &    0.000054406 21942 90990 91465  &  20.58988 67749 18085 64961  \\
  \end{tabular}
   \caption{The sum of $1/n$ where $n$ has \textit{exactly}, or \textit{at most}, two of each digit, base $b$}
   \label{Ta:TwoOfEachDigit}
  \end{center}
\end{table}

\textbf{Example 6.}
Consider the sum of $1/n$ where $n$ has exactly \textit{three} of each decimal digit.
This sum is about 0.000010064 61009 49479 35986 66614 43569 .
In order to obtain this result, the algorithm also computed the $4^{10}$ sums of $1/n$ where $n$ has zero, one, two, or three occurrences of each digit.
This computation took about a week on a 2011-era PC.
The sum of $1/n$ where $n$ has \textit{at most} three of each digit, also part of \verb+iSum+'s output, is about 34.00782 14342 15657 46392 34516 96316 .




Table \ref{Ta:ThreeOfEachDigit} shows the results of similar calculations for bases 2 through 10.

\begin{table}[ht]
 \begin{center}
  \begin{tabular}{ r l l }
  $b$ &     \textit{exactly} three of each digit  & \textit{at most} three of each digit \\ \hline
   2  &    0.23033 73608 64593 58053    &  3.94782 64352 57692 74143  \\
   3  &    0.094468 85247 32610 28516  &  7.56978 83546 45729 25186  \\
	 4  &    0.030515 93397 43385 68330  &  11.29252 15268 42711 55967  \\
	 5  &    0.0088407 89873 17666 97690  &  15.06421 24935 92985 51104  \\
	 6  &    0.0024115 45892 41655 21677  &  18.85575 07210 82141 71749  \\
	 7  &    0.00063327 05112 20590 45079  &  22.65231 88430 89428 99556  \\
	 8  &    0.00016201 57262 45651 16005  &  26.44590 01958 47338 12836  \\
	 9  &    0.000040672 01839 59913 75471  &  30.23194 87285 48128 30492  \\
	10  &    0.000010064 61009 49479 35987  &  34.00782 14342 15657 46392  \\
  \end{tabular}
   \caption{The sum of $1/n$ where $n$ has \textit{exactly}, or \textit{at most}, three of each digit, base $b$}
   \label{Ta:ThreeOfEachDigit}
  \end{center}
\end{table}


\textbf{Example 7.}
In textbooks, most examples of convergent series have sums that are less than about 3.
Frequently-encountered subseries of the harmonic series, such as $\sum 1/2^k$, $\sum 1/k^2$, and $\sum 1/k!$, come to mind.
However, the examples above suggest how to construct nontrivial subseries of the harmonic series that have arbitrarily large, but computable, sums.
For example, let's compute the sum of $1/n$ where $n$ has three 0's: \verb+iSum[0, 3]+ gives this printed output:
\begin{verbatim}
    sum = 23.025851037148538
     sum with exactly 0 occurrences of 0 = 23.103447909420542
     sum with exactly 1 occurrence  of 0 = 23.026735341569127
     sum with exactly 2 occurrences of 0 = 23.025860682735520
     sum with exactly 3 occurrences of 0 = 23.025851037148538
     sum with at most 3 occurrences of 0 = 92.181894970873727
\end{verbatim}
Each of the four values of ``sum with exactly $k$ occurrences'' is, very roughly, $10 \ln(10)$ (see Section \ref{S:S10}).
The last line is the total of the four lines directly above it, so this total is about $40\ln(10)$.
In general, the sum of $1/n$ where $n$ has \textit{at most} $M$ 0's would be roughly $(M + 1) \cdot 10 \ln(10)$.
We can make this sum exceed 1000 by taking $M \ge \frac{1000}{10 \ln(10)} - 1 \approx 43$.
In fact, \verb+iSum[0, 43]+ shows that the sum with \emph{at most} 43's is about $1013.21593$.


For $M = 434$, the sum exceeds 10000: the sum with at most 434 0's is about 10016.32364 57764 01861 09739.
The Kempner series considered in \cite{Schmelzer} can also have arbitrarily large sums.
For example, the sum of $1/n$ where $n$ has no occurrence of the digit string ``314159'' is about 2302582.33386 37826 07892 02376.

\textbf{Example 8.}
Table \ref{Ta:ZeroOneTwo} displays, for each decimal digit $d$, the sum of $1/n$ where $n$ has zero, one, or two occurrences of $d$.
The table also shows, for example, that the sum of $1/n$ where $n$ has \textit{at most} two 9's
is about $22.92067 \ldots + 23.04428 \ldots + 23.02604 \ldots  \approx 68.99100 \ldots$ .

\begin{table}[ht]
 \begin{center}
  \begin{tabular}{ r l l l }
  $d$ &     zero occurrences  & one occurrence  & two occurrences \\ \hline

0  &  23.10344 79094 20541 61603  &  23.02673 53415 69126 96109  &  23.02586 06827 35519 97642 \\
1  &  16.17696 95281 23444 26658  &  23.16401 85942 72832 04085  &  23.02727 62863 56005 71224 \\
2  &  19.25735 65328 08072 22453  &  23.08826 06627 56342 39334  &  23.02648 59737 68470 65598 \\
3  &  20.56987 79509 61230 37108  &  23.06741 08819 30230 10242  &  23.02627 31906 67935 05960 \\
4  &  21.32746 57995 90036 68664  &  23.05799 24133 81824 39576  &  23.02617 78853 92600 17317 \\
5  &  21.83460 08122 96918 16341  &  23.05272 88945 30117 49904  &  23.02612 48753 15647 60861 \\
6  &  22.20559 81595 56091 88417  &  23.04940 99732 95500 55704  &  23.02609 15498 64887 12587 \\
7  &  22.49347 53117 05945 39818  &  23.04714 61901 98641 85083  &  23.02606 88649 14415 07436 \\
8  &  22.72636 54026 79370 60283  &  23.04551 39079 82155 53342  &  23.02605 25308 45693 67648 \\
9  &  22.92067 66192 64150 34816  &  23.04428 70807 47848 31968  &  23.02604 02659 61243 78845 \\

  \end{tabular}
   \caption{The sum of $1/n$ where $n$ has exactly zero, one, or two of each digit, $d$}  
   \label{Ta:ZeroOneTwo}
  \end{center}
\end{table}

For $d = 1, 3, 4, 5, 6$, and 7, the values in Table \ref{Ta:ZeroOneTwo} for zero occurrences of $d$ are greater by 1
in the last decimal place than the corresponding values in \cite{Baillie} because the values in \cite{Baillie} were \emph{truncated}.
The values in Table \ref{Ta:ZeroOneTwo} have been calculated to much higher accuracy and are correctly rounded to 20 decimals.

\textbf{Example 9.}
Here, we consider a base greater than 10, namely, base 16 (``hexadecimal'').
Although programmers often use A through F to represent the hexadecimal digits whose (decimal) values are 10 through 15, \emph{Mathematica} just uses the values 10 through 15.

The sum of $1/n$ where $n$ has \emph{one} digit in base 16 whose (decimal) value is 15, can be calculated with \verb+iSum[15, 1, -1, 16]+.
The result is $44.368334748702945$.
Note that, with $b = 16$, this is close to $b \ln(b) \approx 44.361419555836400$;
see Section \ref{S:S10} for more about this approximation.

At what point does the sum reach, say, 40?
\verb+iPartialSumThreshold[15, 1, 40, -1, 16]+ returns $\{60, 39.803395224593713, 61, 40.031792484814502\}$.
The base is 16, so the results mean that the sum for $n < 16^{60} \approx 1.76685 \cdot 10^{72}$
is about 39.80340, while the sum for $n < 16^{61} \approx 2.82696 \cdot 10^{73}$ is about 40.03179.

The partial sum over $n < 16^6$ is \verb+iPartialSum[15, 1, 6, -1, 16]+ $= 2.385144904046052$.
To check this value, run \verb+iPartialSumCheck[15, 1, 6, -1, 16, 1]+.
This takes far longer than \verb+iPartialSum[ ]+, but the result is the same to 14 decimal places.

If the \emph{sixth} parameter is set to 1, \verb+iPartialSumCheck[ ]+ also prints the values of $n < 16^2$ that have one digit 15 in base 16.
These values (in base 10) are: \{15, 31, 47, 63, 79, 95, 111, 127, 143, 159, 175, 191, 207, 223, 239, 240, 241, 242, 243, 244, 245, 246, 247, 248, 249, 250, 251, 252, 253, 254\}.
For example, $127_{10} = 7 \cdot 16 + 15$, so the two digits of $127_{10}$ in base 16 are 7 and 15.
You can check that these are, indeed, the correct values of $n$ to use in the sum.

\textbf{Example 10.}
Consider the series that has no 9's in base 1000.
\verb+iSum[9, 0, -1, 1000]+ returns a sum of 6802.410165253090787.
The first ``missing'' denominator is 9.
In base 10, the next two missing denominator are 19 and 29.
In base 1000, the next two missing denominators are $19_{1000} = 1 \cdot 1000 + 9 = 1009$
and $29_{1000} = 2 \cdot 1000 + 9 = 2009$.

\section{Details of partial sums} \label{S:PartialSumDetails}




An Irwin series consists of a sum of $1/n$ where we limit the number of occurrences of one or more digits in $n$.

When calculating such a sum, we might wonder about the following:

Consider the $n$ in such a series that have, say, $k$ digits, that is, where $10^{k-1} \leq n < 10^k$.
Do these $n$ constitute a smaller and smaller proportion of \emph{all} integers with $k$ digits?
Because these series are convergent, we generally expect that these restricted partial sums will decrease as $k$ increases.
But do they always decrease?
Or, do these partial sums increase for a few values of $k$, \emph{then} decrease?
If the latter, what value of $k$ maximizes the partial sum over $n$ having $k$ digits?

First, here is a trick you can use to obtain data for the partial sums over $n$ that have $k$ digits:
As you increase the ``print level'' in the \verb+irwinSums.m+ package, \verb+iSum[ ]+ prints out more and more intermediate data.
If you set the level to 4 using \verb+setPrintLevel[4]+, then you will see partial sums for \emph{every} value of $k$
that is used in the calculation until convergence occurs.

\textbf{The ``no 9'' series.}
For this series, we will calculate
(a) the \emph{number} of terms in the series that have exactly $k$ digits;
(b) the \emph{proportion} of numbers with $k$ digits that are in the series, relative to the total number of integers having $k$ digits;
and (c) the partial sum over just those $n$ that have exactly $k$ digits, as a function of $k$.

(a) From Equation \eqref{E:Occurrences0}, we know that there are $8 \cdot 9^{k - 1}$ numbers with $k$ digits having no 9's.
This number always \emph{increases} as $k$ increases.

(b) There are $10^k - 10^{k - 1} = 9 \cdot 10^{k - 1}$ numbers with exactly $k$ digits.
Therefore, the \emph{proportion} of terms with $k$ digits that are in the ``no 9'' series is
\[
\frac{8 \cdot 9^{k - 1}}{9 \cdot 10^{k - 1}} = \frac{8}{9} \cdot \left( \frac{9}{10} \right)^{k-1} \, .
\]
This always \emph{decreases} as $k$ increases.

(c) The partial sum of the series over $n < 10^k$ is given by \verb+iPartialSum[9, 0, k]+.
The difference
\begin{verbatim}
  iPartialSum[9, 0, k] - iPartialSum[9, 0, k - 1]
\end{verbatim}
gives the sum of the series over $10^{k-1} \leq n < 10^k$, that is, those $n$ that have exactly $k$ digits.
These partial sums always \emph{decrease} as $k$ increases.

\textbf{The ``two 9'' series.}
Now consider the sum of $1/n$ where $n$ has two 9's.
We can examine the numbers (a), (b), and (c) that we looked at in the ``no 9'' series.
Here, things are a little different.

(a) The \emph{number} of terms in the series that have exactly $k$ digits:
Equation \eqref{E:OccurrencesNotZeroM} shows that, with $m = 2$, there are
\[
 \binom{k - 1}{1} \cdot 9^{k - 2} + 8 \cdot \binom{k - 1}{2} \cdot 9^{k - 3} =
(k - 1) \cdot 9^{k - 2} + 4 \cdot (k - 1)(k - 2) \cdot 9^{k - 3}
\]
numbers with $k$ digits that have two 9's.
This number always increases as $k$ increases.

(b) The \emph{proportion} of numbers with $k$ digits that are in the series:
This proportion is
\[
r_k = \frac{  (k - 1) \cdot 9^{k - 2} + 4 \cdot (k - 1)(k - 2) \cdot 9^{k - 3}  }{9 \cdot 10^{k - 1}} \, .
\]
After simplifying, this reduces to
\begin{equation} \label{E:rkTwoNines}
r_k = \frac{ (k - 1) (4 k + 1)}{9 \cdot 10^2} \cdot \left( \frac{9}{10} \right)^{k-3} \, .
\end{equation}

\begin{figure}[ht] 
\centering
\begin{minipage}[t]{\figureMinipageWidth}
\centering
\includegraphics[width=\picDblWidth]{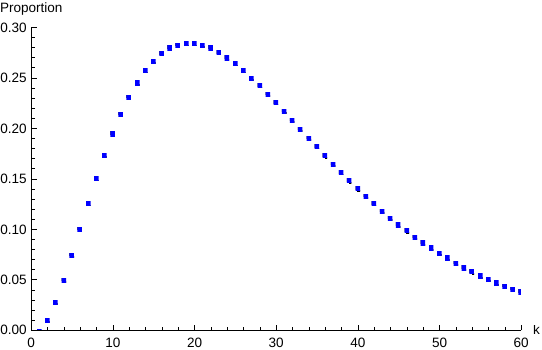}
\caption{$r_k$ from Equation \eqref{E:rkTwoNines}}
\label{fig:rkPlot}
\end{minipage}\hfill
\begin{minipage}[t]{\figureMinipageWidth}
\centering
\includegraphics[width=\picDblWidth]{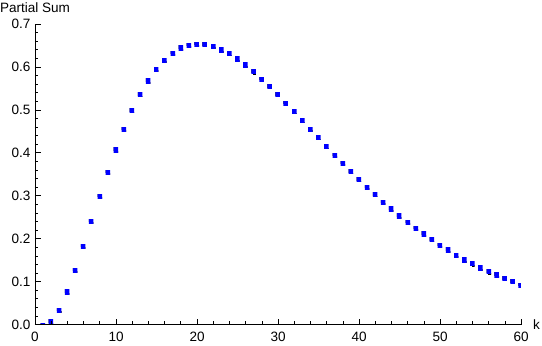}
\caption{Partial sum: $n$ has $k$ digits}
\label{fig:partialSumPlot}
\end{minipage}
\end{figure}

As a function of a \emph{real} variable $k$, $r_k$ has a maximum of about 0.2854799408 at $k = 19.3779991279$.
But here, we are interested only in \emph{integer} values of $k$.
Looking at the numeric values, we see that the largest value of $r_k$ is 0.2853651091 at $k = 19$.
Figure \ref{fig:rkPlot} shows the plot of  $r_k$ for $1 \leq k \leq 60$.

(c) The partial sum restricted to those $n$ that have exactly $k$ digits:
The difference
\begin{verbatim}
  iPartialSum[9, 2, k] - iPartialSum[9, 2, k-1]
\end{verbatim}
gives the sum over just those $n$.
The plot of these values as a function of $k$ is in Figure \ref{fig:partialSumPlot}.
Notice that they \emph{increase} for a while, then decrease.
The maximum partial sum is about 0.6560609728; this is the sum of $1/n$ where $n$ has 20 digits.

By the time we reach $k = 60$, $r_k = 0.0389448094$, and the partial sum over those $n$ that have 60 digits is 0.0952515127 .
The sum of the entire series is \verb+iSum[9, 2]+ $= 23.026040265961244$ .

\textbf{The ``three 9'' series.}
Again, we use Equation \eqref{E:OccurrencesNotZeroM}, but with $m = 3$.
Here, the proportion of $n$ with $k$ digits, relative to the total number of $n$ with $k$ digits, is
\[
r_k = \frac{K_d(k, 3)}{9 \cdot 10^{k-1}} \, .
\]
Among \emph{integer} values of $k$, this $r_k$ reaches a maximum of about 0.2361884380 at $k = 29$.
The maximum partial sum is about 0.5432931624; this is the sum of $1/n$ where $n$ has 30 digits.

By the time we reach $k = 60$, $r_k = 0.0838327871$, and the partial sum over those $n$ that have 60 digits is 0.2019732126 .
The sum of the entire series is \verb+iSum[9, 3]+ $ = 23.025852998372444$ .

\textbf{The ``one 8 and one 9'' series.}
Here, we have one occurrence of each of two non-zero digits.
We use Equation \eqref{E:OccurrencesOneAOneB} to count the number of these $n$ that have $k$ digits.
Again, let $r_k$ for this series be the proportion of these $n$, relative to the total number of integers with $k$ digits.
This $r_k$ reaches a maximum of 0.1514609778 at $k = 9$.
This maximum does not depend on \emph{which} two non-zero digits these are.
The maximum partial sum, 0.3463458110, occurs for $k = 10$ digits.
This sum \emph{does} depend on the two digits being 8 and 9; the sum will be different for other pairs of digits.

By the time we reach $k = 60$, $r_k = 0.0000828040$, and the partial sum over those $n$ that have 60 digits is 0.0002235618 .
The sum of the entire series is \verb+iSum[{8, 9}, {1, 1}]+ $ = 5.756763686415333$ .

\ifthenelse {\boolean{BKMRK}}
  { \section{The \texorpdfstring{$s_{10}$}{s10} series}\label{S:S10} }
  { \section{The s10 series}\label{S:S10} }

Define $s_{m}$ to be the sum of $1/n$ where $n$ has exactly $m$ zeros.

The $s_{10}$ series begins with the rather small terms
\[
s_{10} = \frac{1}{10000000000} + \frac{1}{20000000000} + \frac{1}{30000000000} + \cdots \, .
\]
The sum of this series is about $23.025850929940457$.

A. D. Wadhwa \cite{Wadhwa} showed that as $m$ inreases, $s_m$ strictly decreases, and that $s_m > 19.28$ for all $m$.

With \verb+iSum[0, 0]+, \verb+iSum[0, 1]+, etc., we can calculate that
\begin{align*}
s_0 & \approx 23.10344 \ 79094 \ 20542 \\
s_1 & \approx 23.02673 \ 53415 \ 69127 \\
s_2 & \approx 23.02586 \ 06827 \ 35520 \\
s_3 & \approx 23.02585 \ 10371 \ 48538 \\
s_4 & \approx 23.02585 \ 09311 \ 18602 \\
s_5 & \approx 23.02585 \ 09299 \ 53404 \\
s_{10} & \approx 23.02585 \ 09299 \ 40456 \ 84018 \ 19892 \ 30669 \, .
\end{align*}

This last result was obtained with \verb+iSum[0, 10, 30]+.
Notice that $s_{10} - 10 \ln(10) \approx 2.075 \cdot 10^{-21}$.
Question: Why is this difference so small?
Answer: Bakir Farhi \cite{Farhi} has shown that, for \emph{any} digit $d$, if $s_m(d)$ is the sum of $1/n$
where $n$ has exactly $m$ occurrences of $d$,
then, as $m$ approaches infinity, $s_m(d)$ converges, from above, to
\[
10 \ln(10) \approx 23.02585 \ 09299 \ 40456 \ 84017 \ 99145 \ 46844 \ .
\]

In base $b$, the sums approach $b \ln(b)$.
For example, the sum of $1/n$ where $n$ has ten 0's in base 2 is calculated with \verb+iSum[0, 10, -1, 2]+.
This is about 1.386298363144833, which is about $2 \ln(2) + 4.002 \cdot 10^{-6}$.
With larger bases, the sum with ten 0's is even closer to $b \ln(b)$:
In base 16, the sum is about $16 \ln(16) + 7.724 \cdot 10^{-26}$.

\textbf{When does the partial sum exceed $10 \ln(10)$?}

According to Farhi's result, for any finite $m$, the sum of the series $s_m$ is a little more than $10 \ln(10)$.
For example, $s_{10} \approx  10 \ln(10) + 2 \cdot 10^{-21}$.
Let's calculate at what point the $s_{10}$ series exceeds $10 \ln(10)$.

If we run \verb+iPartialSumThreshold[0, 10, 10*Log[10]]+, we get a warning that the threshold (23.025850929940457) is too close to the sum of the series.
This is because this calculation uses only 15 decimal places (the default).
But the sum differs from $10 \ln(10)$ by about $2.075 \cdot 10^{-21}$, so we will need more than 20 decimals.
Let's use about 30.
We can either use \verb+N[10*Log[10], 30]+ to specify the threshold to 30 digits (28 decimals),
\begin{verbatim}
  x = iPartialSumThreshold[0, 10, N[10*Log[10], 30] ]
\end{verbatim}
or we can use the fourth parameter of \verb+iPartialSumThreshold[ ]+ to specify the number of decimals for the function to use, as with
\begin{verbatim}
  x = iPartialSumThreshold[0, 10, 10*Log[10], 30]
\end{verbatim}
Both of these will return a value of \verb+x+ that is approximately
\[
\{759, 23.025850929940456840179738941229, 760, 23.025850929940456840179937336511\} \, .
\]
This shows that, in order for the partial sum to exceed $10 \ln(10)$, we need at least some denominators having 760 digits.
The two partial sums (the second and fourth components of \verb+x+) are about $10 \ln(10) - 1.76 \cdot 10^{-22}$ and $10 \ln(10) + 2.28 \cdot 10^{-23}$.

The \emph{number of terms} required to reach this point can be found with Equation \eqref{E:OccurrencesZeroM}.
For $m = 10$, that Equation says that number of terms in the $s_{10}$ series which have $k$ digits, is
\begin{equation} \label{E:s10kCountTerms}
9 \cdot \binom{k - 1}{10} \cdot 9^{k - 11} = \binom{k - 1}{10} \cdot 9^{k - 10} \, .
\end{equation}
Summing this for $k \leq 759$ and $k \leq 760$, the two sums are $9.75 \cdot 10^{736}$ and $8.89 \cdot 10^{737}$.

So, there are about $9.75 \cdot 10^{736}$ terms in the series that are less than $10^{759}$;
we need \emph{more} terms than this for the partial sum to exceed $10 \ln(10)$.

There are about $8.89 \cdot 10^{737}$ terms in the series that are less than $10^{760}$;
the sum of these terms exceeds $10 \ln(10)$ by about $2.28 \cdot 10^{-23}$.

Table \ref{Ta:s10} shows partial sums of the $s_{10}$ series at various thresholds.
For each threshold $t$, the return value from \verb+iPartialSumThreshold[0, 10, t, 30]+ is of the form $\{k, p_k, k+1, p_{k+1}\}$.
$p_k$ and $p_{k+1}$ are the partial sums for $n < 10^k$ and $n < 10^{k+1}$, so that $p_k < t < p_{k+1}$.
``Terms'' is the \emph{total} number of terms in the series for $n < 10^{k+1}$.
So, the sum of this many terms equals $p_{k+1}$ and \emph{exceeds} the threshold, $t$.


\begin{table}[ht]
 \begin{center}
  \begin{tabular}{ r r l l l }
  Threshold ($t$)   &  $k$ &    $p_k$  & $p_{k+1}$ & Terms to $10^{k+1}$ \\ \hline
\rule{0pt}{12pt}  
$10^{-5}$ & 19 & $0.90470 \cdot 10^{-5}$ & $1.81536 \cdot 10^{-5}$ & $3.4 \cdot 10^{14}$ \\
$0.1$ & 45 & $0.09660$ & $0.11553$ & $7.9 \cdot 10^{43}$ \\
 1 & 62 & $0.99442$ & $1.09930$  & $4.5 \cdot 10^{61}$ \\ 
 2 & 69 & $1.88538$ & $2.04340$ & $6.8 \cdot 10^{68}$ \\ 
 3 & 75 & $2.95004$ & $3.15420$ & $8.8 \cdot 10^{74}$ \\ 
 4 & 79 & $3.81000$ & $4.04243$ & $1.0 \cdot 10^{79}$ \\ 
 5 & 83 & $4.77800$ & $5.03508$ & $1.1 \cdot 10^{83}$ \\ 
10 & 101 & $9.96759$ & $10.2708$ & $1.3 \cdot 10^{101}$ \\ 
15 & 118 & $14.8350$ & $15.0906$ & $1.1 \cdot 10^{118}$ \\ 
20 & 144 & $19.8954$ & $20.0250$ & $5.7 \cdot 10^{143}$ \\ 
21 & 153 & $20.9130$ & $21.0067$ & $4.1 \cdot 10^{152}$ \\ 
22 & 168 & $21.9911$ & $22.0416$ & $2.2 \cdot 10^{167}$ \\ 
23 & 232 & $22.9998$ & $23.0015$ & $7.2 \cdot 10^{229}$ \\ 
23.02585 & 372 & $23.025849925$ & $23.025850001$ & $3.4 \cdot 10^{365}$ \\ 
$t = 10 \ln(10)$ & 759 & $t - 1.76 \cdot 10^{-22}$ & $t + 2.28 \cdot 10^{-23}$ & $8.9 \cdot 10^{737}$
  \end{tabular}
   \caption{ $s_{10}$: Thresholds, partial sums, and total terms with $n < 10^{k+1}$}
   \label{Ta:s10}
  \end{center}
\end{table}

From Equation \eqref{E:s10kCountTerms}, the \emph{proportion} of $k$-digit integers that are in the series is
\[
\frac{ \binom{k - 1}{10} \cdot 9^{k - 10} } {9 \cdot 10^{k-1} } \, .
\]
This proportion reaches a maximum of about 0.1318653468 at both $k = 100$ and $k = 101$.
The sum over $n$ having exactly 100 digits is about 0.3037008645.
This is larger than the sum over $n$ having any other number of digits.
These statements can be confirmed using the trick mentioned near the beginning of Section \ref{S:PartialSumDetails}.

\section{The googol series} \label{S:Googol}

One googol is $10^{100}$, which can be written as 1 with 100 zeros.
In the sum of $1/n$ where $n$ has exactly 100 zeros, the first term is the very tiny 1/googol.
The first few terms are
\[
s_{100} = \frac{1}{10^{100}} + \frac{1}{2 \cdot 10^{100}}+ \frac{1}{3 \cdot 10^{100}} + \cdots \, .
\]
Although the first term is tiny and later terms are even smaller,
the sum of this series is $10 \ln(10) + 1.0074572170 \ldots \cdot 10^{-197} \approx 23.02585$ .
Table \ref{Ta:s100-d} has a more accurate value of this sum.

\emph{How large must $N$ be} in order to make the partial sum over $n < N$ exceed $10 \ln(10)$?

The function \verb+iPartialSumThreshold[ ]+ gives us the answer.
 
To make the partial sum exceed $10 \ln(10)$, we need at least some denominators $> 10^{7226}$ :
\begin{verbatim}
  x = iPartialSumThreshold[0, 100, 10*Log[10], 220]
\end{verbatim}
The result (after about half an hour on a 2020-era laptop) is, roughly,
\begin{verbatim}
  x = {7225, 23.02585 ..., 7226, 23.02585 ...}
\end{verbatim}
The partial sums (the second and fourth elements of \verb+x+) have more than 200 digits after the decimal point;
both are close to $10 \ln(10)$, and they differ from each other in the $201^{\text{st}}$ decimal place.

Subtracting $10 \ln(10)$ from the second and fourth elements of \verb+x+, we get:
\begin{verbatim}
  x[[2]] - 10 Log[10] = -2.83888407016547054930*10^-200
  x[[4]] - 10 Log[10] = 8.544994808119088314378*10^-199
\end{verbatim}
So, the partial sum over $n < 10^{7225}$ is about $10 \ln(10) - 2.84 \cdot 10^{-200}$,
while the partial sum over $n < 10^{7226}$ is about $10 \ln(10) + 8.54 \cdot 10^{-199}$.

Note that we must do these calculations to more than 200 decimal places in order to verify that the second partial sum really does exceed $10 \ln(10)$.


\emph{How many terms} are in this series with, say, $n < 10^{7226}$?
From Equation \eqref{E:OccurrencesZeroM}, the number of terms in this series that have exactly $k$ digits is
\begin{equation}  \label{E:s100kCount1}
\binom{k - 1}{100} \cdot 9^{k - 100} \, .
\end{equation}
We sum this expression for $k = 101$ through $k = 7226$.
The sum, the number of terms with $n < 10^{7226}$, is about $4 \cdot 10^{7027}$.

The \emph{proportion} of $k$-digit integers that occur as denominators in the $s_{100}$ series is
\[
\frac{ \binom{k - 1}{100} \cdot 9^{k - 100} }{9 \cdot 10^{k-1}} \, .
\]
This proportion reaches a maximum of about 0.04201 67909 at both $k = 1000$ and $k = 1001$.
This makes sense: in a random number with $k = 1000$ digits, we expect that about $k/10 = 100$ of them to be zero;
if $k$ was much larger, we would expect more than 100 digits to be zero.

\emph{How large is the sum} over just those $n$ that have 1000 digits?
We can estimate the partial sum of the harmonic series over $n$ having $k$ digits with the integral
\begin{equation} \label{E:IntegralEstimate}
\sum_{n = 10^{k-1}}^{10^k} \frac{1}{n}
\approx \int_{10^{k-1}}^{10^k} \frac{1}{x} \; dx = \ln(10) \approx 2.3025850930 \, .
\end{equation}
The proportion of 1000-digit $n$ in the harmonic series that are in $s_{100}$ is about 0.04201 67909.
Assuming these terms in $s_{100}$ are randomly scattered among all integers with 1000 digits,
we can estimate the partial sum of the $s_{100}$ series over 1000-digit denominators, as about
$\ln(10) \cdot 0.0420167909 \approx 0.09674 \ 72364$.


We can compute the actual sum over 1000-digit denominators with
\begin{verbatim}
  iPartialSum[0, 100, 1000] - iPartialSum[0, 100, 999]
\end{verbatim}
The result is $0.09674 \ 94386 \ 91394$, not far from the above estimate.
This sum is larger than the sum over $n$ having any \emph{other} number of digits.
These statements can be confirmed using the trick mentioned near the beginning of Section \ref{S:PartialSumDetails}.

It is impossible to directly add up the first googol terms of this series.
But we still might wonder:
\emph{About how large would the sum of the first googol terms be?}

To estimate this, we first use the expression in \eqref{E:s100kCount1}.
A little trial and error shows that
\[
\sum_{k = 1}^{158} \binom{k - 1}{100} \cdot 9^{k - 100} = 7.17 \cdot 10^{98} \, ,
\quad \text{ while } \quad
\sum_{k = 1}^{159} \binom{k - 1}{100} \cdot 9^{k - 100} = 1.76 \cdot 10^{100} \, .
\]
So, in the $s_{100}$ series, one googol terms are reached when $n$ is somewhere between $10^{158}$ and $10^{159}$.

Next, we use \verb+iPartialSum[0, 100, 158]+ and \verb+iPartialSum[0, 100, 159]+
to compute the partial sums for $n < 10^{158}$ and $n < 10^{159}$.
These partial sums are about $3.32943 \cdot 10^{-59}$ and $8.21656 \cdot 10^{-59}$.
Therefore, the sum of the first googol terms would lie somewhere between these two values.

\emph{Thresholds.}
Table \ref{Ta:s100} shows partial sums of this series at various thresholds.
$p_k$ and $p_{k+1}$ are the partial sums for $n < 10^k$ and $n < 10^{k+1}$, so that $p_k < t < p_{k+1}$.
``Terms'' is the \emph{total} number of terms in the series for $n < 10^{k+1}$.
So, the sum of this many terms \emph{exceeds} the threshold, $t$.

The sum of the entire series is $10 \ln(10) + 1.00745 \ldots \cdot 10^{-197}$.
The threshold in last line of Table \ref{Ta:s100} is $t = 10 \ln(10) + 10^{-197}$,
which is (very) slightly less than this sum.

\begin{table}[ht]
 \begin{center}
  \begin{tabular}{ l r l l l }
  Threshold ($t$)   &  $k$ &    $p_k$  & $p_{k+1}$ & Terms to $10^{k+1}$ \\ \hline
\rule{0pt}{12pt}  
$10^{-50}$    &  182 &  $0.66901 \cdot 10^{-50}$   &  $1.34440 \cdot 10^{-50}$ & $2.5 \cdot 10^{132}$ \\
1/10          &  777 &  0.09734         &  0.10103 & $1.6 \cdot 10^{775}$ \\
1/2           &  826 &  0.49087         &  0.50527 & $6.2 \cdot 10^{824}$ \\
 1            &  852 &  0.99153         &  1.01670 & $1.1 \cdot 10^{851}$ \\
 2            &  882 &  1.98058         &  2.02260 & $1.8 \cdot 10^{881}$ \\
10            &  990 &  9.95559         & 10.05188 & $4.2 \cdot 10^{989}$ \\
20            & 1116 & 19.96750         & 20.01555 & $2.1 \cdot 10^{1115}$ \\
22.92068      & 1276 & 22.92027         & 22.92298 & $1.2 \cdot 10^{1274}$ \\
23            & 1327 & 22.99988         & 23.00062 & $3.2 \cdot 10^{1324}$ \\
23.02585      & 1611 & 23.025849977     & 23.025850017 & $1.7 \cdot 10^{1604}$ \\
$t = 10 \ln(10) - 10^{-100}$  & 4692 & $t - 8.12 \cdot 10^{-102}$ & $t + 4.76 \cdot 10^{-103}$ & $3.9 \cdot 10^{4591}$ \\
$t = 10 \ln(10)$  & 7225 & $t - 2.84 \cdot 10^{-200}$ & $t + 8.54 \cdot 10^{-199}$ & $4.0 \cdot 10^{7027}$ \\
$t = 10 \ln(10) + 10^{-197}$  & 7278 & $t - 4.58 \cdot 10^{-201}$ & $t + 2.35 \cdot 10^{-201}$ & $3.1 \cdot 10^{7078}$ \\
  \end{tabular}
   \caption{ $s_{100}$: Thresholds, partial sums, and total terms with $n < 10^{k+1}$}
   \label{Ta:s100}
  \end{center}
\end{table}

How close is the partial sum over say, $n < 10^{10000}$ to the sum of the entire series? Answer: \\
\verb+iSum[0, 100, 330] - iPartialSum[0, 100, 10000, 330]+    
is about $1.6190098640 \cdot 10^{-310}$.
There are about $7 \cdot 10^{9688}$ terms in the $s_{100}$ series with $n < 10^{10000}$.

For each digit $d$, Table \ref{Ta:s100-d} shows the sum of $1/n$ where $n$ has 100 occurrences of $d$.
These were computed using \verb+iSum[d, 100, 250]+.

\begin{table}[ht]
 \begin{center}
  \begin{tabular}{ r l }
  $d$ & Sum \\ 

0 &  $10 \ln{10}$ + 1.00745 72170 67704 21141 82347 02395 53683 90097 45688 70487 $ \cdot 10^{-197}$  \\
1 &  $10 \ln{10}$ + 1.46324 68728 28002 21624 31047 42226 13853 19705 84513 60021 $ \cdot 10^{-195}$  \\
2 &  $10 \ln{10}$ + 6.50985 48660 87751 53297 11600 77857 48564 93107 75257 68892 $ \cdot 10^{-196}$  \\
3 &  $10 \ln{10}$ + 4.32791 33177 18277 79749 26702 46868 66487 36553 86198 09026 $ \cdot 10^{-196}$  \\
4 &  $10 \ln{10}$ + 3.35146 07081 41798 51340 34344 03171 12314 44393 94486 17979 $ \cdot 10^{-196}$  \\
5 &  $10 \ln{10}$ + 2.80859 79174 23800 00239 17516 48643 34458 22747 83369 49527 $ \cdot 10^{-196}$  \\
6 &  $10 \ln{10}$ + 2.46741 99007 53799 94215 23312 26886 65419 13349 17702 22343 $ \cdot 10^{-196}$  \\
7 &  $10 \ln{10}$ + 2.23522 17279 43195 05545 67500 48273 20522 70628 87614 34314 $ \cdot 10^{-196}$  \\
8 &  $10 \ln{10}$ + 2.06805 18229 36195 99843 97677 85425 91950 80493 92989 41553 $ \cdot 10^{-196}$  \\
9 &  $10 \ln{10}$ + 1.94253 90508 24014 22882 12900 35330 10473 39190 27815 69249 $ \cdot 10^{-196}$  \\

  \end{tabular}
   \caption{Sum of $1/n$ where $n$ has exactly 100 occurrences of the digit $d$}
   \label{Ta:s100-d}
  \end{center}
\end{table}

\textbf{Other bases.}
Table \ref{Ta:s100-baseB} displays the sum of $1/n$ where $n$ has 100 0's in bases 2 through 16.

\begin{table}[ht]
 \begin{center}
  \begin{tabular}{ r l }
 $b$ & Sum \\ 
 2 & $2 \ln{2}   + 4.58536 \ 70841 \ 64038 \ 69252 \ 13542 \ 25169 \ 75624 \ 53549 \ 68042 \ 67967 \cdot 10^{-49}$ \\
 3 & $3 \ln{3}   + 6.78035 \ 24742 \ 30102 \ 54437 \ 86740 \ 46140 \ 63559 \ 57239 \ 77853 \ 65479 \cdot 10^{-86}$ \\
 4 & $4 \ln{4}   + 7.35621 \ 92954 \ 32871 \ 32308 \ 23741 \ 02975 \ 37118 \ 70647 \ 16970 \ 28836 \cdot 10^{-113}$ \\
 5 & $5 \ln{5}   + 9.15700 \ 51061 \ 45339 \ 45952 \ 05149 \ 47207 \ 95648 \ 20120 \ 11823 \ 22577 \cdot 10^{-134}$ \\
 6 & $6 \ln{6}   + 9.51961 \ 10946 \ 62454 \ 49834 \ 19110 \ 96561 \ 68278 \ 63718 \ 17705 \ 61758 \cdot 10^{-151}$ \\
 7 & $7 \ln{7}   + 5.09401 \ 50691 \ 25087 \ 27809 \ 09774 \ 86098 \ 40629 \ 53652 \ 62009 \ 41180 \cdot 10^{-165}$ \\
 8 & $8 \ln{8}   + 2.58405 \ 95082 \ 48340 \ 95667 \ 52691 \ 61563 \ 61301 \ 95345 \ 90871 \ 14266 \cdot 10^{-177}$ \\
 9 & $9 \ln{9}   + 4.15887 \ 00874 \ 91506 \ 31089 \ 67547 \ 11372 \ 61989 \ 38749 \ 82215 \ 91809 \cdot 10^{-188}$ \\
10 & $10 \ln{10} + 1.00745 \ 72170 \ 67704 \ 21141 \ 82347 \ 02395 \ 53683 \ 90097 \ 45688 \ 70487 \cdot 10^{-197}$ \\
11 & $11 \ln{11} + 2.16293 \ 58305 \ 26135 \ 30598 \ 14668 \ 39457 \ 89780 \ 05710 \ 20083 \ 19602 \cdot 10^{-206}$ \\
12 & $12 \ln{12} + 2.78884 \ 26002 \ 28004 \ 68628 \ 83230 \ 59167 \ 02605 \ 47370 \ 10712 \ 27928 \cdot 10^{-214}$ \\
13 & $13 \ln{13} + 1.60939 \ 56301 \ 78690 \ 18082 \ 01782 \ 31082 \ 86863 \ 34033 \ 25190 \ 86887 \cdot 10^{-221}$ \\
14 & $14 \ln{14} + 3.31107 \ 11603 \ 84035 \ 49131 \ 48002 \ 87996 \ 71118 \ 18646 \ 41505 \ 08983 \cdot 10^{-228}$ \\
15 & $15 \ln{15} + 2.02945 \ 18537 \ 45221 \ 71388 \ 89799 \ 29120 \ 08194 \ 48027 \ 78886 \ 11775 \cdot 10^{-234}$ \\
16 & $16 \ln{16} + 3.20854 \ 13415 \ 41349 \ 56954 \ 55639 \ 35881 \ 31757 \ 70742 \ 73057 \ 52983 \cdot 10^{-240}$ \\
  \end{tabular}
  \caption{Sum of $1/n$ where $n$ has exactly 100 0's in base $b$}
  \label{Ta:s100-baseB}
 \end{center}
\end{table}

\section{Rational or irrational?} \label{S:Irrat}

\textbf{Are the sums we consider here rational or irrational?}
The series in Examples 4 - 6 in Section \ref{S:Examples} have a finite number of terms, and are therefore rational.
But there is no particular reason to think the other sums are rational.
And, because there are more irrationals than rationals, any bets should be placed on the sums being irrational.

However, we do know this: a theorem of Borwein \cite[Theorem 1]{Borwein} states that if $b$ is an integer such that $|b| > 1$ and $c$ is a non-zero rational number, then
\[
\sum_{n=1}^{\infty} \frac{1}{b^n + c} 
\]
is irrational (assuming that $b^n + c$ is not zero for any $n \ge 1$).
If we take $b = 10$ and $c = -1$, this theorem implies that
\[
\sum_{n=1}^{\infty} \frac{1}{10^n - 1} = \frac{1}{9} + \frac{1}{99} + \frac{1}{999} + \cdots
\]
is irrational.
Multiply this series by $9/d$ and we see that if the denominators of a series consist of any single non-zero digit $d$, then the sum is irrational.
The same is true in other bases $b$.


Let $x_d$ be the sum of $1/n$ where $n$ has one occurrence of $d$.
For each digit $d = 0$, 1, ..., 9, the author computed $x_d$ to 5010 decimal places.
(See Table \ref{Ta:ZeroOneTwo} for these sums to 20 decimal places.)
For none of these sums did the Mathematica function \verb+Rationalize+ find a rational number $p/q$ such that
\[
| x_d - \frac{p}{q}  | < \frac{10^{-4}}{q^2}
\]

For $d = 0$ and $d = 9$, the author also computed the corresponding sums to 10010 decimal places.
\verb+Rationalize+ also could not find a $p/q$ meeting the above criterion for either of these sums.

\textit{Mathematica} also has a function, \verb+RootApproximant+, which takes as input a number $c$ and attempts to find a polynomial with integer coefficients that has $c$ as a root.
If $c$ is a number known to $n$ digits, then one can usually find such a polynomial for which the total number of digits in the coefficients is about $n$, even if $c$ is not algebraic.
However, if the total number of digits in the coefficients is much less than $n$, this strongly suggests that $c$ is, indeed, algebraic, with $c$ being a root of that polynomial.
All ten 5010-decimal place approximations $x_d$ with $d = 0$, 1, ..., 9 were tested with \verb+RootApproximant+.
In no case did \verb+RootApproximant+ find a polynomial of degree $\le 10$ where the total number of digits in the coefficients was less than 5000.
Likewise, for the two 10010 decimal place approximations with $d = 0$ and $d = 9$, \verb+RootApproximant+ could not find polynomials of degree $\le 10$ where
the total number of digits in the coefficients was less than 10000.

\section{A challenge to the reader}\label{S:Unsolved}

The algorithm presented in this paper allows one to sum a series whose denominators contain a chosen number of occurrences of one or more digits.
For example, the sum of $1/n$ where $n$ has exactly one 3 and one 5 is about 5.754336150131606.
But how do we compute the sum of $1/n$ where $n$ has, say, exactly one occurrence of a multi-digit integer such as 35?

For Kempner's series, the algorithm in \cite{Baillie} allows one to sum a series when one digit is excluded from the denominators.
Like the algorithm presented here, that algorithm is based on the way that $i$-digit denominators in the series gave rise to the $(i+1)$-digit denominators.

Schmelzer \cite{Schmelzer} generalized the algorithm in \cite{Baillie} to handle the case where a \emph{multi-digit} integer is excluded from the denominators.
For example, \verb+kSum[35]+ in the \verb+kempnerSums.m+ package shows that the sum of $1/n$ where $n$ has no occurrence of 35 is about 227.909716596753556 .
Schmelzer's algorithm works by partitioning the $i$-digit denominators into several subsets, then extrapolating to $(i+1)$-digit denominators on each subset separately.
Can the algorithm presented here be generalized in a similar way to compute the sum of $1/n$ where $n$ has, say, $m$ occurrences of a multi-digit integer?

Let's consider the counts of integers which have $k$ digits, and which have one occurrence of some 2-digit integer.
Let $d$ and $e$ be digits, with $d \neq 0$ and $e \neq d$.
As the 2-digit integer ranges from 00 through 99, there are four cases:
The two-digit integer can be of the form $00$, $0d$, $dd$, or $de$.

The last column of Table \ref{Ta:TwoDigitCounts} shows the counts of integers having $k$ digits
and which have one occurrence of a two-digit integer of the form $de$.
For example, there are 278 4-digit integers with one occurrence of 35.
But there's nothing special about 35.
The counts are the same for \emph{any} two-digit integer $de$ provided $d \neq 0$ and $d \neq e$.

Likewise, the counts where the integer $0d$ takes any of the nine values $01$ through $09$, are equal.
The counts where the integer $dd$ takes any of the nine values $11, 22, \ldots, 99$, are equal.

\begin{table}[ht]
 \begin{center}
  \begin{tabular}{ r r r r r }
 $k$  &    $00$  &      $0d$  &      $dd$  &       $de$  \\ \hline

 2  &         0  &         0  &         1  &          1  \\
 3  &         9  &         9  &        17  &         19  \\
 4  &       162  &       180  &       243  &        278  \\
 5  &      2349  &      2682  &      3141  &       3642  \\
 6  &     30618  &     35460  &     38394  &      44863  \\
 7  &    376164  &    439227  &    452466  &     531317  \\
 8  &   4448358  &   5221080  &   5197041  &    6122876  \\
 9  &  51221727  &  60326964  &  58567131  &   69156804  \\
10  &            &            &            &  769184205
  \end{tabular}
   \caption{Counts of $k$-digit integers having one $00$, $0d$, $dd$, or $de$.}
   \label{Ta:TwoDigitCounts}
  \end{center}
\end{table}

The sum of $1/n$ over $n < 10^{10}$ where $n$ has exactly one 35, is about 1.05002 02131 69342 49813.
The partial sums where $n$ has 2 through 10 digits, are about:
0.02857, 0.05223, 0.07474, 0.09673, 0.11825, 0.13934, 0.15998, 0.18020, and 0.199985.

There is one 4-digit integer with \emph{two} 35's: namely, 3535.
There are 29 5-digit integers with two 35's: $13535, 23535, \ldots, 93535$;
$35035, 35135, \ldots, 35935$;
finally, $35350, 35351, \ldots, 35359$.
The counts of integers with two 35's, and which have 4 through 10 digits, are: 1, 29, 567, 9283, 137126, 1893294, and 24918690.

The sum of $1/n$ over $n < 10^{10}$ where $n$ has exactly two 35's, is about 0.020469 94051 23425 46036.
The partial sums where $n$ has 4 through 9 digits, are about:
0.00028289, 0.00080452, 0.0015490, 0.0025082, 0.0036752, 0.0050435, and 0.0066066.

These sums were computed by ``brute force'', that is, by summing the appropriate terms with denominators less than $10^9$,
not using the algorithm in Section \ref{S:Algorithm}.

These partial sums over $n$ having $k$ digits are increasing for $k \leq 10$.
This suggests that $1.05002 \ldots$ and $0.020469 \ldots$ are not good approximations to the sums of the entire series.

The previous version of this preprint made an incorrect comparison between the sums of $1/n$ where $n$ has one 35 in bases 10 and 100.

In fact, the numbers in base 100 that have one digit whose value is (decimal) 35,
\emph{do not} constitute a proper subset of numbers that contain one 35 in base 10.

For example, (decimal) 135135 has one base-100 digit equal to 35: the base-100 digits are $\{13, 51, 35\}$.
However, this number in base 10 has two 35's, not one.
So, (decimal) 135135 is in the set of integers with one 35 in base 100, but is \emph{not} in the corresponding set for base 10.

Professor Burnol at the University of Lille in France provided a key idea \cite{Burnol-Kempner}, \cite{Burnol-Irwin} to this author that enables a very rough estimate of the sums discussed here.
The key is to use \emph{generating functions}.

Given enough coefficients, \textit{Mathematica} can compute a generating function having those coefficients.
For example, using seven of the the counts in the last column of Table \ref{Ta:TwoDigitCounts}, we can do this:
\begin{verbatim}
  g = FindGeneratingFunction[{1, 19, 278, 3642, 44863, 531317, 6122876}, x]
\end{verbatim}
The generating function is
\begin{align*}
  g & = \frac{1 - x}{(1 - 10 x + x^2)^2} \\
 & = 1 + 19 x + 278 x^2 + \dots + 6122876 x^6 + 69156804  x^7 + 769184205 x^8 + 8451616295  x^9 + \dots \, .
\end{align*}
This generating function reproduces all nine entries in the last column in Table \ref{Ta:TwoDigitCounts}.
For now, let's assume that this generating function also produces the correct counts for all larger $k$.

This \textit{Mathematica} code computes coefficients of $x^0$ through $x^{5000}$:
\begin{verbatim}
  maxPower = 5000
  ser = Series[g, {x, 0, maxPower}] ;
  coeffs = Table[SeriesCoefficient[ser, k], {k, 0, maxPower}] ;
\end{verbatim}

The $k^{\text{th}}$ coefficient (\verb+coeffs[[k]]+) is the \emph{count} of integers having $k+1$ digits, with one occurrence of the digit string 35.

Therefore, the \emph{proportion} of $k$-digit integers with one 35 is the $(k-1)^{\text{st}}$ coefficient divided by the total number of $k$-digit integers, which is $9 \cdot 10^{k-1}$.
In \textit{Mathematica}, this is
\begin{verbatim}
  prop[k_] := coeffs[[k - 1]]/(9 * 10^(k - 1))
\end{verbatim}

If we examine these proportions, we see that they reach a maximum of about 0.373526 at $k = 99$.
By the time we reach $k = 5000$, the proportion is less than $10^{-20}$.

We can follow the same procedure we used with the Googol series.
See Equation \eqref{E:IntegralEstimate} and the paragraph which follows that Equation.
The sum of $1/n$ where $n$ has $k$ digits, and where $n$ has one digit string 35, is about
\begin{verbatim}
  Log[10] * prop[k]
\end{verbatim}

We can now estimate the sum over all $k$.
Here, we compute each proportion to 20 decimals.
\begin{verbatim}
  one35 = Log[10] * Sum[N[prop[k], 20], {k, 2, 5000}]
\end{verbatim}
This estimate for the sum of the entire series is about 230.25850 92994 04568 40.

If we instead sum only to $k = 4000$, the estimate is 230.25850 92994 04546 28.
Keep in mind that these are only rough estimates anyway.

Note that this is a \emph{generic} estimate based only on the \emph{counts} in Table \ref{Ta:TwoDigitCounts}, not on the particular digit string 35.
The \emph{actual} sums for different digit strings of the form $de$, such as 12, 35, or 98 would, of course depend on the particular digit string.

With the limited data in Table \ref{Ta:TwoDigitCounts}, \textit{Mathematica} cannot get the generating functions for the $00$ and $0d$ cases.
However, we can get a generating function for the $dd$ case where $d \neq 0$:
\begin{verbatim}
  gdd = FindGeneratingFunction[{ 1, 17, 243, 3141, 38394, 452466, 5197041 }, x]
\end{verbatim}
we find that
\begin{align} \label{E:gdd}
  gdd & = \frac{1 - x}{(1 - 9 x - 9 x^2)^2} \\
 & = 1 + 17 x + 243 x^2 + \dots + 5197041 x^6 + 58567131  x^7 + 650385369 x^8 + 7138637001  x^9 + \dots \, .
\end{align}
Our estimate for the $dd$ case is:
\begin{verbatim}
maxPower = 10000
serdd = Series[gdd, {x, 0, maxPower}];
coeffsdd = Table[SeriesCoefficient[serdd, k], {k, 0, maxPower}];
propdd[k_] := coeffsdd[[k-1]]/(9 * 10^(k - 1))
\end{verbatim}
We now get an estimate for the sum of $1/n$ where $n$ has one 11 (or one 22, or one 33, etc.)
\begin{verbatim}
  onedd = Log[10] * Sum[N[propdd[k], 35], {k, 2, maxPower}]
\end{verbatim}
This gives 230.25850929940456840179914546843642.
Why is this sum so close to $100 \cdot \ln(10)$?

Professor Burnol pointed out that the infinite sum is \emph{exactly} $100 \cdot \ln(10)$.
First, from Equation \eqref{E:gdd}, we have that $gdd(1/10) = 9000$.
Second, in the \textit{Mathematica} code, \verb+onedd/Log[10]+ is
\[
  \frac{1}{90} + \frac{17}{900} + \frac{243}{9000} + \frac{3141}{90000} + \dots
= \frac{1}{90} \left( 1 + \frac{17}{10} + \frac{243}{100} + \frac{3141}{1000} + \cdots \right)
= \frac{1}{90} \cdot gdd \left( \frac{1}{10} \right) = 100 \, .
\]
One could also ask: What is the sum of $1/n$ where $n$ has exactly \emph{two} (or three, or four) occurrences of, say, 35?
What is the sum where $n$ has, say, two 00's and one 35?

It would be nice to have an algorithm to accurately sum series like these.
However, the author has not pursued this matter further.

\section{Acknowledgments}
The author would like to thank Professor Jean-Fran\c{c}ois Burnol at the University of Lille for finding an error and for very helpful suggestions.

\bigskip



\noindent\textnormal{Email: bobbaillie@frii.com; State College, PA 16803}




\begin{thebibliography}{99}


\bibitem{Baillie}
  Robert Baillie,
  Sums of reciprocals of integers missing a given digit,
  \textit{American Mathematical Monthly} \textbf{86}, no. 5, (May, 1979), pp. 372--374.\\
  Available through JSTOR at \url{https://dx.doi.org/10.2307/2321096} . \\
  The values in Table 1 in this reference are \emph{truncated} to 20 decimals;
  the values in Table \ref{Ta:ZeroOneTwo} above are \emph{rounded}. \\
  Erratum: For digits 1, 2, 3, and 4, the $10^{27}$ on page 373 should be $10^{24}$.
  See ERRATA, \textbf{86}, no. 10 (December, 1980), p. 866.

 \bibitem{Boas}
  Ralph P. Boas, Jr. and John W. Wrench, Jr.,
  Partial Sums of the Harmonic Series,
  \textit{American Mathematical Monthly} \textbf{78} (October, 1971), pp. 864--870.\\
  Available through JSTOR at \url{https://dx.doi.org/10.2307/2316476} .

\bibitem{Borwein-Borwein}
  Jon Borwein and Peter Borwein,
  Evaluation of Sums of Reciprocals of Fibonacci Sequences, in
  \textit{Pi \& the AGM: A Study in Analytic Number Theory and Computational Complexity} (1987),
  Wiley, New York, pp. 91--101.

\bibitem{Borwein}
  Peter Borwein, On the Irrationality of Certain Series,
  \textit{Math. Proc. Cambridge Philos. Soc.} \textbf{112} (1992), 141--146.\\
  Available for U.S. \$36.00 at \url{https://dx.doi.org/10.1017/S030500410007081X} .

\bibitem{Burnol-Kempner}
  Jean-Fran\c{c}ois Burnol,
  Moments in the exact summation of the curious series of Kempner type,
  \\ \url{https://arxiv.org/abs/2402.08525} .

\bibitem{Burnol-Irwin}
  Jean-Fran\c{c}ois Burnol,
  Moments for the summation of Irwin series,
  \url{https://arxiv.org/abs/2402.09083} .

\bibitem{Farhi}
  Bakir Farhi,
  A Curious Result Related to Kempner's Series,
  \textit{American Mathematical Monthly} \textbf{115} (December, 2008), pp. 933--938.\\
  Available through JSTOR at \url{https://www.jstor.org/stable/27642640} ;\\
  a preprint is at \url{https://arxiv.org/abs/0807.3518} .

\bibitem{Foregger}
  T. Foregger,
  Helping Professor Umbugio. Solution to Problem E2533,
  \textit{American Mathematical Monthly} \textbf{83} (August/September, 1976), pp. 570--571.

\bibitem{HardyAndWright}
  G. H. Hardy and E. M. Wright,
  \textit{An Introduction to the Theory of Numbers},
  Oxford, 4th edition, 1960.

\bibitem{Havil}
  Julian Havil,
  \textit{Gamma: Exploring {E}uler's Constant},
  Princeton University Press, Princeton, NJ 2003.

\bibitem{Irwin}
  Frank Irwin,
  A Curious Convergent Series,
  \textit{American Mathematical Monthly}, \textbf{23} (May, 1916), pp. 149--152.\\
  Available free through JSTOR at \url{https://dx.doi.org/10.2307/2974352} .

\bibitem{Kempner}
  A. J. Kempner,
  A Curious Convergent Series,
  \textit{American Mathematical Monthly}, \textbf{21} (February, 1914), pp. 48--50.\\
  Available free through JSTOR at \url{https://dx.doi.org/10.2307/2972074} .

\bibitem{Mathologer}
  Mathologer,
  700 years of secrets of the Sum of Sums (paradoxical harmonic series),
  (``Chapter 6'', which includes a discussion of Kempner and Irwin series, starts at 29:43). \\
  \url{https://www.youtube.com/watch?v=vQE6-PLcGwU} .

\bibitem{Pondiczery}
  E. S. Pondiczery,
  Elementary Problem E2533,
  \textit{American Mathematical Monthly} \textbf{82} (April, 1975), p. 401.

\bibitem{Schmelzer}
  Thomas Schmelzer and Robert Baillie,
  Summing a Curious, Slowly Convergent Series,
  \textit{American Mathematical Monthly} \textbf{115} (June/July, 2008), pp. 525--540. \\
  Available through JSTOR at \url{https://www.jstor.org/stable/27642532} ;\\
  \textit{Mathematica} code that implements the algorithm in the Schmelzer/Baillie article is at\\
  \url{https://library.wolfram.com/infocenter/MathSource/7166} .

\bibitem{BoasOEIS}
  Sequence A004080,
  The On-Line Encyclopedia of Integer Sequences, OEIS Foundation, \\
  \url{https://oeis.org/A004080} .  

\bibitem{BoasOEIS2}
  Sequence A082912,
  The On-Line Encyclopedia of Integer Sequences, OEIS Foundation, \\
  \url{https://oeis.org/A082912} .  

\bibitem{Wadhwa}
  A. D. Wadhwa,
  Some Convergent Subseries of the Harmonic Series,
  \textit{American Mathematical Monthly} \textbf{85} (October, 1978), pp. 661--663.
  Available through JSTOR at \url{https://dx.doi.org/10.2307/2320338} .

\bibitem{smbc-cartoon}
  Zack Weinersmith,
  Saturday Morning Breakfast Cereal, \\
  \url{https://www.smbc-comics.com/index.php?id=3777} .

\bibitem{Weisstein-Digit}
  Eric Weisstein,
  \textit{Digit}, MathWorld--A Wolfram Web Resource, \\
  \url{https://mathworld.wolfram.com/Digit.html} .

\bibitem{Weisstein-Lambert}
  Eric Weisstein,
  \textit{Lambert Series}, MathWorld--A Wolfram Web Resource, \\
  \url{https://mathworld.wolfram.com/LambertSeries.html} .

\bibitem{NegativeBinomialSeries}
  Eric Weisstein,
  \textit{Negative Binomial Series}, MathWorld--A Wolfram Web Resource, \\
  \url{https://mathworld.wolfram.com/NegativeBinomialSeries.html} .

\bibitem{Wells}
  David Wells,
  \textit{The Penguin Dictionary of Curious and Interesting Numbers},
  Penguin Books, London, 1997.

\bibitem{EB-Constant}
  Erdős–Borwein constant,
  \textit{Wikipedia}, \\
  \url{https://en.wikipedia.org/wiki/Erd%C5%91s%E2%80%93Borwein_constant} .



\end{thebibliography}
\end{document}